\documentclass{article}
 
\usepackage{graphicx} 
\usepackage{amsmath}
\usepackage{amssymb}
\usepackage{enumitem}
\usepackage{algorithm2e}
\usepackage[utf8]{inputenc}
\usepackage{booktabs, multirow}
\usepackage{comment}
\usepackage{natbib}
 \bibpunct[, ]{(}{)}{,}{a}{}{,}%
\usepackage[left=1in, right=1in, top=1in, bottom=1in]{geometry}
\usepackage{tabularx}
\usepackage{yhmath}
\usepackage{xcolor}
\usepackage{subcaption}

\title{Enhancements of Fragment Based Algorithms for Vehicle Routing Problems}
\date{}

\author{
  Lucas Sippel\\
  School of Mathematics and Physics, The University of Queensland\\
  \texttt{l.sippel@uq.edu.au}
  \and
  Michael Forbes\\
  School of Mathematics and Physics, The University of Queensland\\
  \texttt{m.forbes@uq.edu.au}
}

\newcommand{\w}{\omega}
\newcommand{\W}{\Omega}
\newcommand{\f}{\mathtt{f}}
\newcommand{\F}{\mathtt{F}}
\newcommand{\FF}{\mathcal{F}}
\newcommand{\BB}{\mathcal{B}}
\newcommand{\N}{\mathtt{N}}
\newcommand{\V}{\mathtt{V}}
\newcommand{\aaa}{\mathtt{a}}
\newcommand{\A}{\mathtt{A}}
\newcommand{\vv}{\mathtt{v}}

\newcommand{\rr}{\mathtt{r}}
\newcommand{\R}{\mathtt{R}}

\renewcommand{\P}{\mathtt{P}}
\newcommand{\p}{\mathtt{p}}
\newcommand{\bp}{\bar{\mathtt{p}}}
\renewcommand{\l}{\mathtt{l}}

\newcommand{\bL}{\bar{L}}
\newcommand{\bl}{\bar{\l}}

\newcommand{\bff}{\bar{f}}
\newcommand{\HH}{\mathtt{H}}

\newcommand{\bc}{\bar{c}}

\newtheorem{theorem}{Theorem}

\newtheorem{definition}{Definition}

\begin{document}

\maketitle

\begin{abstract}
The method of fragments was recently proposed, and its effectiveness has been empirically shown for three specialised pickup and delivery problems. We propose an enhanced fragment algorithm that for the first time, effectively solves the Pickup and Delivery Problem with Time Windows.  Additionally, we describe the approach in general terms to exemplify its theoretical applicability to vehicle routing problems without pickup and delivery requirements.  We then apply it to the Truck-Based Drone Delivery Routing Problem Problem with Time Windows.  The algorithm uses a fragment formulation rather than a route one.  The definition of a fragment is problem specific, but generally, they can be thought of as enumerable segments of routes with a particular structure.  A resource expanded network is constructed from the fragments and is iteratively updated via dynamic discretization discovery.  Additionally, we introduce two new concepts called \emph{formulation leveraging} and \emph{column enumeration for row elimination} that are crucial for solving difficult problems.  These use the strong linear relaxation of the route formulation to strengthen the fragment formulation.  We test our algorithm on instances of the Pickup and Delivery Problem with Time Windows and the Truck-Based Drone Delivery Routing Problem with Time Windows. Our approach is competitive with, or outperforms the state-of-the-art algorithm for both.
\end{abstract}

\section{Introduction}\label{introduction}
Vehicle routing problems (VRPs) aim to calculate a set of vehicle routes that service geographically distributed requests at minimum cost.  They were first introduced by \cite{firstvrp} to determine how to route fuel trucks from a depot to service demand at fuel stations.  

One of the simplest VRPs is the Capacitated Vehicle Routing Problem (CVRP).  It can be modelled on a complete directed graph $G = (V, E)$ where $V = V'\cup\{0,n+1\}$ is the set of vertices and $E$ is the set of edges.  Each vertex $i\in V'$ represents a request with weight $q_i$.  Vertices $0$ and $n+1$ represent the depot at the start and the end of vehicle routes.  Each edge $(i,j)\in E$ has cost $c_{ij}$.  An infinite fleet of homogeneous vehicles is assumed, each with capacity $Q$.  In the CVRP, a route is an elementary path $\rr=(i_1=0,i_2,...,i_k=n+1)$ such that the vehicle capacity is not violated, $\sum_{j=1}^k q_{i_j}\leq Q$.  The cost of a route is the sum of its edge costs.  The CVRP seeks to minimize the cost of the routes used to service each request $i\in V'$ exactly once.

There are numerous extensions to this that model real world applications.  For example, request service time windows can be included in which case the problem is called the Vehicle Routing Problem with Time Windows (VRPTW).  The Pickup and Delivery Problem with Time Windows (PDPTW) further extends the VRPTW by coupling requests into pickup and delivery pairs.Further qualifications are added to the definition of a route with each extension.  Thus, we define a route in the following generic way.
\begin{definition}
    A \emph{route} is an elementary path $\rr = (i_1=0,i_2,...,i_k=n+1)$ in $G$ that satisfies the problem specific rules of the VRP in question.  
\end{definition}
Further extensions such as multiple depots and finite heterogeneous vehicle fleets have also been considered, though these require a more complex definition of a route.  The methods we proposed in this work assume the above definition for simplicity, though they can be extended to solve other VRPs.  The reader is referred to \cite{vrptax1,vrptax2, vrptax3} and \cite{vrptax4} for extensive summaries of VRP variants.

VRPs can be solved using \emph{arc-flow} formulations where binary variables $x_e$ represent a vehicle traversing edge $e\in E$.  These have a polynomial number of variables and constraints and a weak linear relaxation.  \cite{dwdecomp} give a decomposition method which can be used to transform an arc-flow formulation to a set partitioning formulation \citep{spp} which includes a variable for every route.  This formulation consequently has an exponential number of variables, but usually a stronger linear relaxation than the arc-flow formulation.

\subsection{The Route Formulation}
The set partitioning, or route formulation has a binary variable $\lambda_\rr$ for each route $\rr$ in the set of all possible routes $\R$.  For each request $i\in V'$, we let $a^i_\rr$ be one if route $\rr$ services request $i$ and zero otherwise. Thus, the formulation, SPP, is the following:
\begin{align}
    \min& & &\sum_{\rr\in\R} c_\rr\lambda_\rr & &\label{rbfo}\\
    \text{s.t.}& & &\sum_{\rr\in\R}a^i_\rr\lambda_\rr = 1, & &\forall i\in V'\label{rbfc1}\\
    & & &\lambda_\rr\in\{0,1\}, & &\forall \rr\in\R\label{rbfc2}.
\end{align}
Objective function (\ref{rbfo}) minimizes the cost of the routes used.  Constraints (\ref{rbfc1}) ensure each request is covered once.  The domain of the variables is defined by constraints (\ref{rbfc2}).  Whilst many vehicle routing problems require extra constraints (for example, a limit on the number of vehicles), we omit them for the sake of clarity in the rest of Sections \ref{introduction} and \ref{algorithm}.

As stated by \cite{bnpvrps}, an advantage of formulating VRPs in this way is that many of the constraints are implicitly satisfied by each route and need not be included in the formulation.  This usually results in route formulations having better lower bounds than their arc-flow counterparts.  The price paid is a huge number of variables which are often impossible to include all at once.  Branch-and-price-and-cut (BPC) algorithms are designed to circumvent this.

\subsection{Column Generation for BPC Algorithms}
In this section we cover the main component of BPC algorithms, namely, column generation (CG).  We refer the reader to \cite{bnp1}, \cite{cg}, \cite{bnp2} and \cite{bnpvrps} for in depth reviews of BPC algorithms.

A BPC algorithm is a branch-and-bound (BB) algorithm where lower bounds are calculated via CG at nodes of the search tree.  At each node, the linear relaxation of SPP superimposed with the node's set of branching constraints is referred to as the master problem (MP).  We refer to the linear relaxation of SPP before any branching decisions are made as the root MP.  CG solves MP by iteratively solving a restricted MP (RMP) and then a pricing subproblem.  At the root node the RMP associated with SPP and the MP is given by
\begin{align}
    \min& & &\sum_{\rr\in\R'} c_\rr\lambda_\rr & &\label{rmpo}\\
    \text{s.t.}& & &\sum_{\rr\in\R'}a^i_\rr\lambda_\rr = 1, & &\forall i\in V'\label{rmpc1}\\
    & & &\lambda_\rr\geq 0, & &\forall \rr\in\R'\label{rmpc2},
\end{align}
where $\R'$ is a small subset of the routes in $\R$.  

Finding an optimal primal solution, $(\lambda_\rr)_{\rr\in\R'}$,  yields an optimal dual solution $(\pi_i)_{i\in V'}$ to RMP's dual linear program (LP).  The feasible MP solution that is obtained by setting $\lambda_\rr=0$ for all $\rr\in\R\backslash\R'$ is optimal if all routes $\rr\in\R\backslash\R'$ have non-negative reduced cost, $\bar{c}_\rr = c_\rr - \sum_{i\in V'}a^i_\rr\pi_i$.  Hence, the pricing problem (PP) is solved to find negative reduced cost routes to add to $\R'$:
\begin{equation}
    \min_{\rr\in\R\backslash\R'} c_\rr - \sum_{i\in V'}a^i_\rr\pi_i.
\end{equation}
If solving PP yields one or multiple negative reduced cost routes, they are added to $\R'$ and the process repeats.  Otherwise, it terminates with an optimal MP solution.  

Subproblem PP is an elementary shortest path problem with resource constraints (ESPPRC) on $G$.  The resources are quantities, such as time, load or distance, that vary as a route is traversed according to resource extension functions (REFs) which are defined for each resource and for each edge $(i,j)\in E$.  Each resource has lower and upper bounds that must be satisfied.  If the resource constraints permit cycles, one must carefully disallow them or relax the elementarity requirement, in which case the PP is solved as a shortest path problem with resource constraints (SPPRC).  The SPPRC can be solved by a pseudo-polynomial algorithm whilst the ESPPRC is NP-hard.  However, \cite{ESPPRC1} show that using non-elementary routes may result in significantly lower quality bounds from MP. Thus, significant focus has been placed on designing elementarity relaxations that balance the computational complexity of PP and the quality of MP's lower bounds. See \cite{bnpvrps} for details.  Note that when SPPRC is used, the definition of $a^i_\rr$ is modified to be the number of times route $\rr\in\R$ visits request $i\in V'$.  

An SPPRC is solved by maintaining $\P$, a set of \emph{resource feasible} paths in G that each start at $0$.  For efficiency, the set of paths is represented by a tree structure of labels.  The path represented by a label can be recovered from the tree structure and so we use label $\l$ and path $\p$ interchangeably.  A path's label $\l$ stores $i_\l$ its head vertex, $\bc_\l$ its reduced cost, and a vector $L(\l) = (r_\l)_{r\in R}$ of the vehicle's current resource consumption values, where $R$ is the set of resources in the SPPRC.  At each request $i$ in $V$, the feasible set for resource $r$ is $U_i^r$.  We assume that $U_i^r$ is an interval $[\alpha_i^r,\beta_i^r]$ for all numeric resources though this can be extended to a finite union of disjoint intervals. A label $\l$ is resource feasible if $r_\l\in U_{i_{\l}}^r$ for all $r\in R$.  The algorithm proceeds by removing a label $\l$ from $\P$ and extending it from its head vertex along edges $(i_\l,j)\in E$.  Let $f_{ij}^r$ denote the REF for edge $(i,j)\in E$ and resource $r\in R$.  Thus, the new label after extending is $(j,\bc_\l + c_{i_\l j} + \pi_{i_\l/2} + \pi_j/2,(f_{i_\l j}^r(r_\l))_{r\in R})$ and the new label is only stored in $\P$ if it is resource feasible.  

As in \cite{vrpsolver}, we say a numerical resource $r$ is \emph{disposable} if the consumption can be increased at will without the cost of the path being affected.  In this case, the REF satisfies
\begin{equation}
    f_{ij}^r(r_L) = \max\{\hat{f}_{ij}^r(r_L), \alpha_j^r\}
\end{equation}
for some real-valued, non-decreasing function $\hat{f}_{ij}^r$.  We refer to all other resources as \emph{non-disposable}.  For use later, we denote the vector of non-disposable and disposable resource consumptions for each forward path $\p$ by $L^1(\p)$ and $L^2(\p)$ respectively.

Efficiency depends on domination.  With disposable REFs a label $\l_1$ dominates another with the same head vertex, $\l_2$, if $L(\l_1)\leq L(\l_2)$ component-wise and $\bc_{\l_1}\leq\bc_{\l_2}$, in which case the label $\l_2$ can be ignored.  Modifications to this dominance rule are needed when non-disposable resources or path structural constraints are present, including when solving the PP as an ESPPRC.  

In many VRPs paths can also be extended backwards starting from the end depot, provided backward REFs can be formulated \citep{refsinv}.  A backward path is a path $\p$ in $G$ that ends at $n+1$ and can be completed feasibly by a vehicle.  We represent backward labels with $\bl$ and the backward REFs by $\bff_{ij}^r$.  When a backward label is extended from vertex $j$ to $i$, the resource consumptions of new backward label $\bl'$ are $\bL(\bl') = (\bff_{ij}^r(r_{\bl}))_{r\in R}$. We denote the vector of non-disposable and disposable resource consumptions for backward path $\bp$ by $\bL^1(\bp)$ and $\bL^2(\bp)$ respectively.

Backward labelling is important in BPC for a CG acceleration technique known as bidirectional labelling \citep{bicg}.  For details, see \cite{bnpvrps}.  It is also used for a well known variable fixing technique \citep{vf2} which we take advantage of in our approach as discussed in Section \ref{formlev}.

Branching decisions become necessary if the root MP solution is not integer.  In this case we cannot take advantage of efficient commercial BB solvers because CG dynamically adds variables to MP.  Instead a bespoke BB algorithm must be implemented which in turn leads to a range of design considerations and enhancements that strongly influence the efficacy of the algorithm (see \cite{bnpvrps} for details).  Such considerations make the implementation error-prone and time consuming.  While advanced BPC algorithms have been implemented for many VRPs, \cite{vrpsolver} state that ``designing and coding each of those complex and sophisticated BPCs has been a highly demanding task, measured on several work-months of a skilled team."  We name some of the considerations and enhancements here to illustrate why achieving a successful BPC implementation is hard: pricing heuristics \citep{pricingheur}, dual variable stabilisation \citep{dualsmoothing}, both robust and non-robust valid inequalities \citep{ssr}, ng-path relaxation and decremental state space relaxation \citep{originaldssr,ngpath}, completion bounds \citep{originalcompbnd}, strong branching and other branching strategies \citep{pricingheur}, and primal heuristics \citep{divingheur}. One main idea underpinning our proposed algorithm is leveraging the high quality lower bound of MP while simultaneously avoiding much of the complexity of implementing BPC algorithms. 

\subsection{Fragment Methods}
Recently, solution methods involving sub-paths of routes have been implemented for specific vehicle routing problems.  Such methods rely on enumeration of the sub-paths, which is possible because they are shorter than routes.  The resulting formulations based on these objects typically have weaker linear relaxations than their route counterparts.  On the other hand, a priori enumeration removes the need for a labelling algorithm at BB nodes and allows the use of commercial solvers.  Thus, the success of such approaches depends on whether this trade-off is worth it.  

\cite{pstep2} describe a family of formulations for the CVRP based on $p$-steps, which are sub-paths of routes up to length $p$.  \cite{frag1} describe a method for the Pickup and Delivery Problem with Time Windows and First-In-Last-Out loading (PDPTWL) based on fragments.  A fragment is a sub-path of a route that starts with a vehicle empty at a pickup location and ends at the first delivery where the vehicle becomes empty.  More recent work by \cite{frag2} and \cite{frag3} improve on fragment algorithms by using shorter sub-paths than fragments when they cannot be enumerated, and introducing valid inequalities to tighten their formulations.  In \cite{frag3} a sophisticated filtering scheme based on reduced cost is employed to reduce the number of variables in the fragment formulation.

In \cite{frag1}, \cite{frag2} and \cite{frag3} a network is constructed where arcs correspond to fragments and nodes are locations with (possibly) extra information attached.  Care is taken to ensure that the network is relaxed.  That is, any route can be represented by a path (or chain of fragments) through the network.  A network flow IP is defined on the network and is solved via branch-and-cut (BC).  Cuts must be lazily added to the formulation, since chains of fragments in integer solutions are not guaranteed to correspond to routes.    

Despite the fact that the approaches of \cite{frag1}, \cite{frag2} and \cite{frag3} are effective on their respective focus problems, a fragment approach that takes the most appropriate features of each still fails to effectively solve the PDPTW.  The fragments of \cite{frag1} and \cite{frag3} are too long to be enumerated for the PDPTW benchmark instances.  The fragments of \cite{frag2} on the other hand are short enough for enumeration.  Indeed they claim that their approach ``is almost directly applicable to the PDPTW” however, this is not the case for numerous reasons.  Firstly, the time windows of some of the PDPTW benchmarks are wide enough to degrade performance of \cite{frag2}'s approach.  Performance degradation as time windows are widened is a highlighted issue in their paper.  Even if this is addressed by employing a time-expanded network as in \cite{frag1} and \cite{frag3}, the method is ineffective.  For hard instances, convergence is hindered by either the need for too many feasibility cuts, or a weak linear relaxation of the fragment formulation, which makes variable fixing and BB ineffective.

The improvements of our method are therefore motivated by the shortcomings of these previous works to efficiently solve the PDPTW.  However, whilst developing the enhancements, it became apparent that the method may be applicable to VRPs that are not pickup and delivery problems.  This would enable solving such VRPs without implementing a bespoke BB scheme as in BPC.  We therefore describe the framework in general terms and use the PDPTW as an illustrative example.

\subsection{Contributions}\label{contrib}
Our contributions to the VRP literature are three-fold:
\begin{enumerate}
    \item We improve upon the fragment methods of \cite{frag1,frag2} and \cite{frag3} with three key enhancements; \cite{ddd}'s dynamic discretisation discovery (DDD), variable fixing with respect to the route formulation as in \cite{vf2}, and a new technique described in \ref{cere}.  While DDD is a technique proposed previously, we extend it to handle more complex networks for the fragment approach.  Even though \cite{vf2}'s variable fixing technique is well known, the observation that it is relevant to fragment approaches is not trivial.  We test our approach on the PDPTW benchmark instances and show its superiority by comparing its performance to a baseline algorithm that combines the most appropriate aspects of \cite{frag1,frag2} and \cite{frag3}.  Our approach can also outperform the state-of-the-art BPC solver for the PDPTW. 
    \item We describe our approach and its enhancements in general terms to exemplify its theoretical applicability to a wide range of VRPs.
    \item Finally, we apply our approach to the Truck-Based Drone Delivery Routing Problem with Time Windows (TDDRP) proposed by \cite{tddrp}, which shows that or approach can be effective for VRPs other than pickup and delivery problems.  We compare our approach's performance to the BPC algorithm for the TDDRP, showing its superiority.
\end{enumerate}

The rest of the paper is structured as follows.  Section \ref{pdptw} introduces the PDPTW and describes some details about how routes are generated in BPC algorithms for it.  This will aid in the description of our method and applying it to the PDPTW.   Section \ref{algorithm} describes our algorithm and uses the PDPTW as an illustrative example to solidify the reader's understanding of each concept.  Section \ref{pdptw_old} presents and analyses computational results of our approach on the PDPTW.  Section \ref{tddrp} describes how our framework can be extended to the TDDRP and tests its performance on benchmark instances.  Finally, Section \ref{conclusion} gives concluding remarks.

\section{The PDPTW}\label{pdptw}
In this section, we introduce the PDPTW which will be used as an illustrative example throughout the description of our overall approach.  We have $G=(V,E)$ where $V=\{0, 2n+1\}\cup P\cup D= \{0,...,2n,2n+1\}$.  Vertices $0$ and $2n+1$ are the start and end depots.  Vertices $1,...,n\in P$ are pickup requests.  Each pickup $p$ has corresponding delivery $p+n\in D$.  The routes must satisfy pairing and precedence constraints that these imply;  the delivery $p+n\in D$ of each pickup $p\in P$ must occur after $p$ on the same route.  This complexity is handled by a non-disposable resource in the labelling algorithm of BPC for the PDPTW as discussed in the next subsection.  Each $i \in V$ has time window $[\alpha_i, \beta_i]$.  The travel time for edge $(i,j)\in E$ is denoted $t_{ij}$.  Pickup $p$ has weight $q_p > 0$. For each pickup $p\in P$, $q_{p+n} = -q_p$ to signify the request being delivered.  There are assumed to be infinite homogeneous vehicles with capacity $Q$.  Convention for the PDPTW is to first minimize the number of routes, then minimize the total travel time.

\subsection{Labels}
We briefly describe the forward labelling algorithm for calculating negative reduced cost PDPTW routes without elementarity.  For the sake of brevity, we do not describe backward labelling.  See \cite{pdptw2} for details. The algorithm has labels of the form $\l = (i_\l, \bc_\l, (t_\l, q_\l, O_\l))$, where $i_\l$ and $\bc_\l$ are as previously described, $t_\l$ is the time, $q_\l$ is the current vehicle load, and $O_\l$ is the set of open delivery requests.  Set $O_\l$ can be represented by binary vector $(O^d_\l)_{d\in D}$ where $O^d_\l$ represents whether or not delivery $d$ is onboard.  The initial label is $(0,0,(0,0,(0)_{d\in D}))$.  We define the REFs for resources $t$, $q$ and each $O^d$. Thus for any extension $j\in V$,
    $f_{ij}^t(t_\l) = \max\{t_\l+t_{ij}, \alpha_{j}\}$,
    $f_{ij}^q(q_\l) = q_\l + q_j$,
and 
\begin{equation}
    f_{ij}^{O^d}(O^d_\l) = \begin{cases}
        O^d_\l - 1 & \text{if j=d},\\
        O^d_\l + 1 & \text{if j=d-n},\\
        O^d_\l &\text{otherwise},
    \end{cases}
\end{equation}

for each $d\in D$.  Clearly, the resource feasibility bounds for $t$ and $q$ at the new location are $U_j^t = [\alpha_j,\beta_j]$, $U_j^q = [0,Q]$.  If the new location is not $2n+1$, $U_j^{O^d} = [0,1]$ for each $d\in D$ maintains the pairing and precedence constraints.  Only labels with no onboard requests ($O^d_\l=0$ for all $d\in D$) can be extended to $2n+1$.  Therefore, $U_{2n+1}^{O^d} = [0,0]$. Note that resource $q$ is redundant because load feasibility can be confirmed by checking $\sum_{d\in D} -q_dO^d_\l \leq Q$.  It is only included for efficiency.  Resource $t$ is disposable whilst resources $O^d$ are not.

\section{The Proposed Algorithm}\label{algorithm}
In this section we describe our enhanced fragment framework for VRPs and use the following small PDPTW instance with three pickup and delivery pairs to illustrate.  Table \ref{pdptw_instance} gives the request weights, time windows and travel times between vertices for the instance.  Vehicles are assumed to have a capacity of 15.  The optimal solution for this PDPTW is thus $(0, 1, 4, 2, 3, 5, 6, 7)$.

\begin{table}[!htp]\centering
\caption{Data for a small PDPTW instance.}\label{pdptw_instance}
\scriptsize
\begin{tabular}{lrrrrrrrrrrrr}\toprule
Request &$q_i$ &$\alpha_i$ &$\beta_i$ &$0$ &$1$ &$2$ &$3$ &$4$ &$5$ &$6$ &$7$ \\\midrule
0 &0 &0 &200 &0 &25.73 &17.37 &25.69 &16.1 &26.72 &10.55 &0 \\
1 &8 &50 &70 & &0 &38.68 &43.65 &41.29 &26.35 &26.35 &25.73 \\
2 &7 &80 &110 & & &0 &9.15 &9.97 &44.02 &27.62 &17.37 \\
3 &8 &100 &130 & & & &0 &18.22 &52.39 &35.36 &25.69 \\
4 &-8 &80 &110 & & & & &0 &40.26 &26.38 &16.1 \\
5 &-7 &160 &190 & & & & & &0 &17.66 &26.72 \\
6 &-8 &170 &200 & & & & & & &0 &10.55 \\
7 &0 &0 &200 & & & & & & & &0 \\
\bottomrule
\end{tabular}
\end{table}

\subsection{The Fragment Approach}\label{ren}
Fragment methods rely on enumerating fragments, which are route sub-paths with (possibly) extra information.
\subsubsection{Fragments for the PDPTW}\label{pdptwfrag}
we combine the idea of \emph{restricted fragments} from \cite{frag2} with the idea of \emph{extended fragments} from \cite{frag3}.  These papers also distinguish fragments from \emph{start arcs} but we consider both of these as fragments.  Therefore, given a PDPTW route $\rr=(i_1,...,i_k)$ a PDPTW fragment $\f$ has path $\p_\f = (i_l,...,i_m)$ for some $1\leq l<m\leq k$ which satisfies one of following:
\begin{enumerate}
    \item Path $\p_\f = (0,p)$ for some pickup $p\in P$.
    \item Path $\p_\f$ starts at a pickup, ends at a pickup or the end depot, and includes one edge from a pickup to a delivery (that is, $i_j\in P$ and $i_{j+1}\in D$ for exactly one $j\in\{l,...,m-1\}$).
\end{enumerate}
It is possible for the vehicle to begin traversing sub-paths of the second type having requests already onboard for delivery.  Fragments must therefore also store the set of onboard deliveries that the vehicle starts with, which we denote $L^1_\f=L^1((i_1,...,i_l))$. 

The small example has 13 such fragments: \{((0,1),$\emptyset$),  ((1,4,2),$\{4\}$), ((1,4,3),$\{4\}$), ((1,4,7),$\{4\}$), ((1,2,4,5,7),$\{4\}$), ((1,2,4,3),$\{4\}$), ((0,2),$\emptyset$), ((2,5,7),$\{5\}$), ((2,3,5,6,7),$\{5\}$), ((0,3),$\emptyset$), ((3,6,7),$\{6\}$), ((3,2,5,6,7),$\{6\}$), ((3,5,6,7),\{5,6\})\}.  The enumeration procedure is a simplified version of that in \cite{frag2}.  Any route can be made from a chain of these fragments as each can be partitioned into sub-paths that satisfy the above rules. 
See \cite{frag2} for a more details.  

Observe that fragment ((3,5,6,7),\{5,6\}) begins with delivery $5$ onboard.  The only fragment that can precede this one in a route is ((1,2,4,3),$\{4\}$) since the vehicle picks up $2$ during the traversal of its path without completing the corresponding delivery $5$.  In the same way, the only fragment that can proceed ((1,2,4,3),$\{4\}$) in a route is ((3,5,6,7),\{5,6\}).  This illustrates why each fragment must contain a starting set of onboard requests; so that chains of fragments obey the pairing and precedence constraints of pickups and deliveries.  

The main idea to understand is that PDPTW fragments are sub-paths that satisfy particular rules, together with a starting set of onboard requests, which is a non-disposable resource of the PDPTW.  The exact starting value of the non-disposable resource must be included for each fragment to ensure that fragment chains represent paths that satisfy the pairing and precedence constraints. 

\subsubsection{Fragments in General} \label{genfrag}For a general VRP, the definition of fragments therefore relies on the subset of its resources that are non-disposable. We denote this subset $R_1\subseteq R$.
\begin{definition}
    Given route $\rr=(i_1,...,i_k)$, a \emph{fragment}, $\f$, consists of two parts:
    \begin{enumerate}
        \item A sub-path $\p_\f = (i_l,...,i_m)$ with $1\leq l < m\leq k$.  
        \item The non-disposable resource vector $L^1_\f = L^1((i_1,...,i_l))$ corresponding to path $(i_1,...,i_l)$.
    \end{enumerate}    
    Also, let $\f^+ = i_l$ and $\f^- = i_m$.  We say $\f$ \emph{covers} requests $i_l,...,i_{m-1}$.  
    
    Notice that in the PDPTW, we do not enumerate all fragments but only those that satisfy rules of our choice which we call \emph{fragment rules}.  These ensure that at least one \emph{appropriate sequence} of fragments exists for each route $\rr=(i_1,...,i_k)\in\R$.  That is, there exists indices $1=j_1<...<j_\kappa=k$  such that $((i_{j_d},...i_{j_{d+1}}), L^1((i_1,...,i_{j_{d}})))$ satisifies the fragment rules for each $d=1,...,\kappa-1$.  Whilst not necessary, it is preferable that the fragment rules permit at most one appropriate sequence for each route.  This avoids symmetry in the fragment IP defined later.  The set of fragments that satisfies the fragment rules is denoted $\F$.
\end{definition} 

\subsubsection{Fragment Resource Extension Functions}
For each resource $r\in R$ and each fragment $\f\in \F$, we define \emph{fragment REF} $\tau_\f^r = f^r_{i_{k-1}i_k}\circ ... \circ f^r_{i_1i_2}$. These calculate the end resource consumption values of the vehicle given that it traverses $\p_\f$.  For non-disposable resources $r\in R_1$, we consider $\tau_\f^r$ a constant rather than a function, because fragments' non-disposable consumptions are fixed.

Let the set of disposable resources be $R_2 = R\setminus R_1$.  Each fragment $\f\in\F$ has a feasible set of starting disposable resource vectors $T_\f$.  That is, given $\f\in\F$, $T_\f$ contains all vectors $L=(L_r)_{r\in R_2}$ that encode feasible starting disposable resource consumption values for the vehicle to traverse $\p_\f$.  More concretely, if $\p_\f=(i_1,...,i_k)$, $T_\f$ contains all vectors $L=(L_r)_{r\in R_2}$ such that $L_r\in U_{i_1}^r$ for all $r\in R_2$ and
    $f^r_{i_{j-1}i_j}\circ ... \circ f^r_{i_1i_2}(L_r)\in U_{i_j}^r$
for $r\in R_2$ and $j=2,...,k$.  

Again we use the small PDPTW instance to illustrate.  The only disposable resource for the PDPTW is time and so $R_2 = \{t\}$.  Using fragment $\f=((3,5,6,7),\{5,6\})$ as an example, observe that $\tau^t_\f(t) = \max\{\max\{\max\{t+52.39, 160\}+17.66, 170\}+10.55,0\}$ and $T_\f = [100,119.4]$, because the vehicle can only pickup $3$ at or after 100 and departure times that are later than 119.4 result in arrival times at end depot $7$ that are later than 200.  Here, 119.4 is calculated as $\min\{\min\{\min\{\beta_7-t_{6,7},\beta_6\}-t_{5,6},\beta_5\}-t_{3,5},\beta_3\} = \min\{\min\{\min\{200-10.55,200\}-17.66,190\}-52.39,130\}$.

\subsubsection{Nodes}
Fragment enumeration yields the set of \emph{nodes} $\HH$.  A node $\eta\in \HH$ is a request $i_\eta\in V$ together with a non-disposable resource vector $L^1_\eta$ such that there exists a fragment $\f$ with $\f^+ = i_\eta$ and $L^1_\f = L^1_\eta$ or $\f^- = i_\eta$ and $(\tau_\f^r)_{r\in R_1} = L^1_\eta$.  In the former case, we say $\f$ starts at $\eta$.  In the latter case we say $\f$ ends at $\eta$.

In the PDPTW, each node $\eta\in\HH$ is a depot or pickup vertex $i_\eta\in P\cup\{0,2n+1\}$ together with a set of onboard deliveries, $L^1_\eta\subset D$.  We illustrate using the small example.  Each of the 13 fragments starts and ends at a node in $\HH = \{(0,\emptyset), (1,\{4\}), (2,\{5\}), (3,\{6\}), (3,\{5,6\}), (7,\emptyset)\}$.

\subsubsection{The Resource Expanded Network}
As previously discussed, the paths that are represented when fragments are connected cannot violate route constraints that are maintained by non-disposable resources because each fragment includes the starting non-disposable resource values in its definition.  However, fragments can still break constraints associated with disposable resources.  For instance, consider fragments $\f_1 = ((1,4,3), \{4\})$ and $\f_2 = ((3,2,5,6,7),\{6\})$ from the small example.  We have $\tau_{\f_1}^t(t) = \max\{\max\{t+41.29,80\}+18.22,100\}$ and so a lower bound on a vehicle's arrival time at pickup 3 given it completes this fragment is $\max\{\max\{50+41.29,80\}+18.22,100\} = 109.51$.  Unfortunately, the feasible starting times for a vehicle completing $\f_2$ fall in set $T_{\f_2} = [100,100.85]$ since later departures than 100.85 lead the vehicle to arrive too late at pickup 2.  Thus connecting these fragments leads to a path that is not resource feasible. 

This motivates the construction of a network that encodes the approximate value of disposable resource consumptions for vehicles as they complete fragments.  In particular, the PDPTW necessitates a network  that encodes time.  Such time expanded networks have been employed by previous fragment approaches for pickup and delivery problems \citep{frag1,frag3}, and in other contexts such as service network design \citep{ddd} and machine scheduling \citep{cgexform}.  Their recent adoption has been specifically catalysed by \cite{ddd}'s observation that they can be applied to problems beyond just those with discrete time points.  We will now give a general description of the network employed in our fragment approach.  We call it the \emph{resource expanded network} because as will become evident throughout the following description, it is theoretically possible for the network to not only encode time, but multiple disposable resources simultaneously.

Given fragments, their REFs, and nodes as previously defined, the resource expanded network $\N=(\V,\W,\A)$ consists of a set of \emph{resourced fragments} $\W$, \emph{resourced nodes} $\V$, and holdover arcs $\A$.
\begin{definition}
    A \emph{resourced fragment} $\w$ is a fragment $\f_\w$ together with a starting disposable resource vector $L_\w^2\in T_\f$.  To simplify notation, we also denote the non-disposable resource vector and the path of its fragment by $L^1_\w = L^1_{\f_\w}$ and $\p_\w = \p_{\f_\w}$ respectively.  Because all starting resource consumptions of resourced fragments are fixed, their true end resource consumptions are known values which we denote $\tau_\w^r$.  These are equal to $\tau_{\f_\w}^r$ and $\tau_{\f_\w}^r(L^2_{\w r})$ for non-disposable and disposable resources $r$ respectively. Finally, we denote the true end non-disposable and disposable resource consumption vectors of $\w$ by $\tau^1_\w$ and $\tau^2_\w$ respectively.  For instance, a resourced fragment associated with fragment $\f_2$ from above has the form $\w = (\f,L^2_\w)$ where $L^2_\w\in[100,100.85]$, $L^1_\w = \{6\}$, and $\f_\w = ((3,2,5,6,7),\{6\})$.  Its true end onboard delivery set is $\tau_\w^O = \tau^1_\w = \emptyset$ and supposing that $L^2_\w=100$, its true end arrival time is $\tau_\w^t =\tau^2_\w= \max\{\max\{\max\{\max\{100+9.15, 80\}+44.02,160\}+17.66,170\}+10.55,0\} = 188.21$.
\end{definition}
\begin{definition}
    A \emph{resourced node} $\vv$ is a node $\eta_\vv\in \HH$ together with a disposable resource vector $L^2_\vv$.  To simplify notation, we let the request and the non-disposable resource vector of its node be $i_\vv = i_{\eta_\vv}$ and $L^1_\vv=L^1_{\eta_\vv}$ respectively.  
    
    Given two resourced nodes $\vv_1$ and $\vv_2$, we say $\vv_1$ is \emph{earlier} than $\vv_2$ if $i_{\vv_1} = i_{\vv_2}$, $L^1_{\vv_1} = L^1_{\vv_2}$, $L^2_{\vv_1}\leq L^2_{\vv_2}$ and there exists at least one resource $r\in R_2$ where the inequality is strict.  We say $\vv_1$ is \emph{later} than $\vv_2$ if $i_{\vv_1} = i_{\vv_2}$, $L^1_{\vv_1} = L^1_{\vv_2}$, $L^2_{\vv_1}\geq L^2_{\vv_2}$ and there exists at least one resource $r\in R_2$ where the inequality is strict.  Let $\mathtt{earlier}(\vv)$ and $\mathtt{later}(\vv)$ be the set of earlier and later resourced nodes than $\vv$ respectively.  Again we solidify understanding of resourced nodes with the small PDPTW example.  Imagine the network includes three resourced nodes (or timed nodes) associated with node $\eta=(3,\{5,6\})$, $\vv_1=(\eta,t_1)$, $\vv_2=(\eta,t_2)$ and $\vv_3=(\eta,t_3)$, and $\alpha_3 \leq t_1 < t_2 < t_3 \leq \beta_3$.  Then $\vv_1$ is earlier than $\vv_2$ whilst $\vv_3$ is later.
    
    Finally, let $\V'$ be the set of resourced nodes $\vv$ with request vertices $i_\vv\in V'$.  We assume $\N$ always has a single resourced node for the start depot $\vv_0$ and end depot $\vv_{n+1}$ whose resource consumption vectors correspond to the trivial path $(0,)$.
\end{definition}
\begin{definition}
    A \emph{holdover arc} leaves one resourced node and arrives at a later resourced node.  In the PDPTW, these connect consecutive resourced (in this case, timed) copies of each node $\eta\in\HH$.
\end{definition}

Each resourced fragment $\w\in\W$ leaves a resourced copy of its fragment's first request $\w^+ = ((\f_\w^+, L^1_\w), L^2_\w)$, and arrives at a resourced copy of its fragment's last request $\w^- = ((\f_\w^-, L'^1_\w), L'^2_\w)$. We say $\w$ starts at $\w^+$ and ends at $\w^-$.  The resource consumption values of these resourced nodes must satisfy specific requirements to ensure every route can be represented by a path through $\N$.  The following properties outline these requirements and inform the network's construction.
\begin{enumerate}[label=Property \arabic*,itemindent=*]
\item For each node there must exist a corresponding resourced node with minimum disposable resource vector.  That is, for each $\eta\in \HH$ there exists $\vv\in\V$ with $\eta_\vv=\eta$ such that $L^2_\vv = (\alpha_i^r)_{r\in R_2}$.  In the small PDPTW instance this means that timed node $\vv=((1,\{4\}), 50)$ must exist for example.
\label{prop1}
\item Given a resourced node $\vv$, a resourced copy of every fragment that can be feasibly completed starting with $\vv$'s resource consumption values must exist.  In other words, for every resourced node $\vv\in\V$ and fragment $\f\in\F$ such that $\f^+ = i_\vv$, $L^1_\f = L^1_\vv$ and $L^2_\vv\in T_\f$, there exists resourced fragment $\w\in\W$ with $\f_\w = \f$ and $\w^+ = \vv$.  Given resourced node $((1,\{4\}),50)$ must exist in the network for the small PDPTW instance, this property enforces that resourced fragment $(((1,4,3),\{4\}),50)$ exists for example.  \label{prop3}
\item Each resourced fragment $\w\in\W$ ends at a resourced node with disposable resource consumption component wise at most the true end consumption of $\w$.  That is, $L'^1_\w=\tau^1_{\w}$ and $L'^2_\w\leq\tau^2_{\w}$.  Take timed fragment $\w = (((1,4,3),\{4\}), 50)$ from the small PDPTW instance for example.  It has true end onboard deliveries $\tau_\w^O = \{6\}$ and arrival time $\tau_\w^t = 109.51$.  It would therefore end at a resourced node of the form $\w^-=((3,\{6\}), t)$ where $t\in[100,109.51]$.  Such an $\w^-$ must exist by \ref{prop1}.
\label{prop2}
\item Each resourced fragments ends at a resourced node for which there is no later alternative.  That is, for every resourced fragment $\w\in\W$, there does not exist resourced node $\vv\in\mathtt{later}(\w^-)$ with $L^2_\vv\leq \tau_\w^2$.  Continuing with resourced fragment $\w = (((1,4,3),\{4\}), 50)$ as an example, suppose that additionally to $\vv=((3,\{6\}),100)$, one other resourced copy $\vv'=((3,\{6\}),105)$ exists. Then $\w$ would end at $\vv'$ rather than $\vv$.\label{prop4}
\item For each pair of resourced nodes $\vv_1$ and $\vv_2$ where $\vv_2$ is later than $\vv_1$, there exists a path of holdover arcs between them. This is satisfied for the PDPTW by connecting consecutive resourced copies of each node.\label{prop5}
\end{enumerate}
We intentionally write these properties similarly to \cite{ddd} to illustrate that resource expanded networks extend the idea of their partially time expanded networks in two ways.  Firstly, resource expanded networks expand over (possibly) multiple disposable resources simultaneously instead of just time.  Additionally, instead of timed copies of edges from $G$ as in \cite{ddd}, we have resourced copies of fragments which are based on sub-paths of routes.  As will become evident in the rest of Section \ref{algorithm}, this second distinction is important for the performance of the fragment approach.  We use the terminology of \cite{ddd} and say that if a resource expanded network satisfies Properties 2 and 3, then it has the \emph{early arrival property}.  Since there are many resource expanded networks that satisfy Properties 1-3, Property 4 is additionally imposed.  It is known as the \emph{longest arc property}.  Observe that when there is more than one disposable resource ($|R_2| > 1$) a lexicographic rounding order must be assigned to the disposable resources to ensure the uniqueness of $\N$.  \ref{prop5} ensures that for each resourced node, all later resourced nodes can be reached via holdover arcs.

\subsubsection{The Fragment Formulation} 
We define the following IP on the resource expanded network once it is constructed.  Let $c_\w$ be the cost of traversing resourced fragment $\w$.  Let $a_\w^i$ be one if $\f_\w$ covers $i$ and zero otherwise.  Let $\W_\vv^+$ and $\W_\vv^-$ denote the set of resourced fragments leaving from and arriving at resourced node $\vv$ respectively and $\A^+_\vv$ and $\A^-_\vv$ denote the same for holdover arcs.  The fragment formulation which we call RFN is a mixed-integer network flow problem on the resource expanded network.  It has binary variable $x_\w$ for each resourced fragment $\w\in\W$ and continuous variable $y_\aaa$ for each holdover arc $\aaa\in\A$.  
\begin{align}
    \min &&&\sum_{\w\in\W} c_\w x_\w& \label{rfno}\\
    \text{s.t.}&&&\sum_{\w\in\W}a_\w^i x_\w = 1 &\forall i\in V',\label{rfnc1}\\
    &&\sum_{\aaa\in\A^+_\vv}y_\aaa + &\sum_{\w\in\W^+_\vv}x_\w = \sum_{\aaa\in\A^-_\vv}y_\aaa + \sum_{\w\in\W^-_{\vv}}x_\w&\forall \vv\in\V',\label{rfnc2}\\
    &&&x_\w\in\{0,1\}&\forall \w\in\W,\label{rfnc3}\\
    &&&y_\aaa\geq 0&\forall \aaa\in\A.\label{rfnc4}
\end{align}
Objective (\ref{rfno}) minimizes the cost of resourced fragments.  Notice that complex costs which are dependent on fragment combinations or whole routes, like duration, are not modelled by (\ref{rfno}).  Constraints (\ref{rfnc1}) ensure each request is covered by exactly one resourced fragment.  Constraints (\ref{rfnc2}) ensure conservation of flow at resourced nodes.  Constraints (\ref{rfnc3}) and (\ref{rfnc4}) define the domain.

Integer solutions to RFN yield \emph{chains} of resourced fragments.  These are essentially sequences of resourced fragments that connect at resourced nodes in the network.
\begin{definition}
    A \emph{chain} is a sequence of resourced fragments $(\w_1,...,\w_k)$ satisfying $\f_{\w_j}^-=\f_{\w_{j+1}}^+$,$ 
        L'^1(\w_j)=L^1(\w_{j+1})$, and
        $L'^2(\w_j)\leq L^2(\w_{j+1})$,
    for $j=1,...,k-1$.  
    
Recall that the resourced fragments have corresponding paths $\p_{\w_j} = (i_1^j,...,i_{m_j}^j)$ for $j=1,...,k$.  Therefore, the chain is a \emph{representation} of the path
    $(i_{1}^1,...,i_{m_1-1}^1)\oplus(i_{1}^2,...,i_{m_2-1}^2)\oplus...\oplus(i_{1}^k,...,i_{m_k}^k)$,
where $\oplus$ is the path concatenation operator.  For simplicity, let $\p_{{\w_1}}\oplus...\oplus \p_{{\w_k}}$ denote $(i_{1}^1,...,i_{m_1-1}^1)\oplus(i_{1}^2,...,i_{m_2-1}^2)\oplus...\oplus(i_{1}^k,...,i_{m_k}^k)$.  We say two chains $(\w^1_1,...,\w^1_{k_1})$ and $(\w^2_1,...,\w^2_{k_2})$ are \emph{equivalent} if $\p_{{\w^1_1}}\oplus...\oplus \p_{{\w^1_{k_1}}} = \p_{{\w^2_1}}\oplus...\oplus \p_{{\w^2_{k_2}}}$.
\end{definition}

Chains can therefore represent routes.  In fact, because of the early arrival property (rounding down end disposable resource consumptions), every route can be represented by a chain, as highlighted in the result below.  This means that RFN is a relaxation of the VRP.
\begin{theorem}\label{repex}
    Suppose we have route $\rr=(i_1,...,i_k)\in\R$.  Then $\rr$ has a representation in $\N$.
\end{theorem}

\begin{figure}[h]
\begin{subfigure}{0.5\textwidth}
    \includegraphics[width=\textwidth]{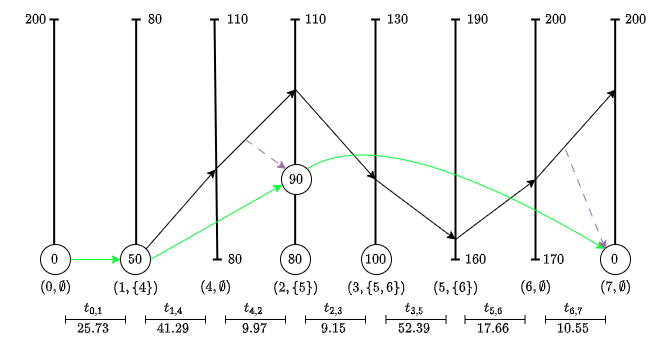}
    \caption{A chain representing a route.}
\end{subfigure}
\hfill
\begin{subfigure}{0.5\textwidth}
    \includegraphics[width=\textwidth]{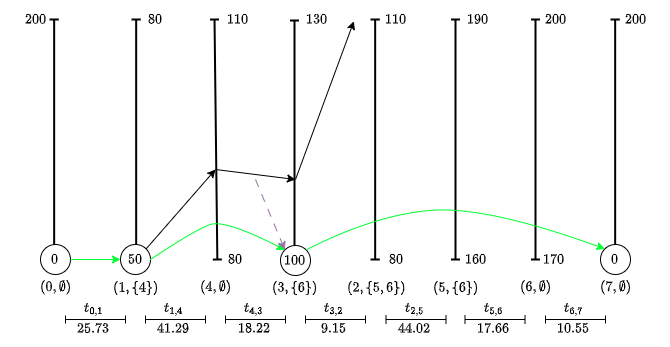}
    \caption{A chain representing an infeasible path.}
\end{subfigure}
\caption{Two RFN solutions for the small PDPTW instance.}
\label{smallpdptw}
\end{figure}

We finish this section by using the optimal route for the small PDPTW instance as an example of a representation and a solution to RFN.  We assume that the network includes the minimal set of resourced nodes as defined by Property 1, plus a second resourced copy of node $(2,\{5\})$ with time $90$.  The left hand side of Figure \ref{smallpdptw} depicts the only representation of the optimal solution for the small PDPTW instance (however it is not the optimal RFN solution as will be discussed in the next section).  The current location and set of onboard deliveries for the vehicle are listed at the bottom of the diagram.  For each, the service time window at the location is depicted.  The circles represent the current resourced nodes in the network.  The solid black arrows represent the actual arrival times of the vehicle at each location and the solid green lines represent the timed fragments of the representation.  The dashed grey lines depict the rounding of time to the latest resourced node at the end of each timed fragment.  Notice that the second timed fragment $\w=(((1,4,2), \{4\}), 50)$ has true arrival time $\tau_\w^t = 101.26$ which is rounded down such that $\w$ ends at the latest appropriate resourced node, $((2,\{5\}),90)$.  The reader can calculate that fragment $\f = ((2,3,5,6), \{5\})$ has feasible set of start times $T_\f = [80,110.25]$ and so by Property 2, a resourced copy leaves the same resourced node.  Observe that because the vehicles actual arrival time at $2$ is also in $T_\f$, the chain represents a resource feasible path.   

\subsubsection{The Relaxation-Size Trade-off}
The rounding down of disposable resource consumptions, however, enable chains to represent paths that are not resource feasible.  We call these \emph{underestimating chains}. An example for the small PDPTW instance is shown on the right hand side of Figure \ref{smallpdptw} made from timed fragments $\w_1=(((1,4,3),\{4\}),50)$ and $\w_2=(((3,2,5,6), \{6\}),100)$.  The true arrival time of the first is $\tau_{\w_1}^t = 109.51$ which is rounded down such that $\w_1$ ends at $((3,\{6\}), 100)$.  Timed fragment $\w_2$ exists because the latest time fragment $((3,2,5,6),\{6\})$ can start is 100.85.  Therefore, the chain exists even though the represented path is not resource feasible. 

Current state-of-the-art fragment methods handle this in two steps:
\begin{enumerate}
    \item They include timed copies of nodes that are a fixed duration $\delta$ apart.  Smaller values of $\delta$ increase the number of timed nodes which removes underestimating chains due to the longest arc property.  For example, if the network depicted in Figure \ref{smallpdptw} included a timed copy of $(3,\{6\})$ at time $105$ then the underestimating chain depicted would be eliminated.
    \item Remaining underestimating chains are eliminated by feasibility cuts in a BC framework.  See \cite{frag3} for details.
\end{enumerate}  
Increasing the number of timed nodes clearly strengthens RFN's linear relaxation because it removes fractional underestimating chains from participating in the solution.  Additionally, it reduces the number of underestimating chain eliminations needed for convergence.  However, it also increases the number of flow conservation constraints (\ref{rfnc3}) and timed fragment variables which may cause inefficiency by slowing down the processing of BB nodes.  We call this the \emph{relaxation-size trade-off}.  

Efficient fragment methods therefore rely on finding a small enough fragment formulation that has a strong enough relaxation to be solved quickly.  For the PDPTW, no such fragment formulation has been found until this work.  Sections \ref{ddd}-\ref{cere} detail three enhancements that enable such a fragment formulation.  The first contributes by reducing the number of resourced nodes needed for RFN to be sufficiently strong, by using DDD.  This alone is not enough for some hard instances however, where sufficient strength is impossible no matter how many resourced nodes are included.  One major contributing factor to the weakness is that non-elementary paths can be represented by fractional chains in the solution of RFN's linear relaxation.  Clearly, basing the fragments on longer paths would reduce the impact of this, and may also remove underestimating chains, resulting in further strengthening of RFN.  In fact, with sufficient increase in the path length of fragments, RFN's linear relaxation would produce a bound equal to that of the route formulation.  See \cite{pstep2}'s discussion on this phenomenon in their p-step formulation of the CVRP.  However, naively increasing fragment path length exponentially increases the number of resourced fragment variables and therefore would be likely to only exacerbate the relaxation-size trade-off.  The second and third enhancements therefore aim to lengthen the paths of fragments without significantly increasing the number of variables.  This is achieved through a combination of dynamically concatenating fragments, and using a well known BPC variable fixing technique to remove many of the resourced fragment variables.

\subsection{Enhancement 1: Dynamic Discretization Discovery}\label{ddd}  

The main drawback of the static networks used by previous fragment methods is that they may include many unnecessary resourced nodes.  On the other hand, some nodes might not have a sufficient number of resourced copies, leading to the elimination of too many underestimating chains.  We handle these problems by implementing \cite{ddd}'s DDD. It works by repeatedly solving RFN and adding resourced nodes to eliminate underestimating chains in integer solutions.  The construction of the network is therefore guided by successive solutions.  

Using DDD to eliminate underestimating chains requires multiple IPs to be solved.  However, DDD has two advantages over BC on a static network.  Firstly, $\N$ can start with few resourced nodes.  Since resourced nodes are only added to remove underestimating chains as needed, DDD often converges to optimality before $|V|$ is large.  This is not the case when a static network is used because a prohibitive number of resource nodes may be needed to ensure that few enough underestimating chains need eliminating.  Furthermore, adding resourced nodes often removes more related underestimating chains than the lifted feasibility cuts in \cite{frag3}.

We leave the full description of DDD to the appendices because it is well known, and only explain the main idea here.  In our approach, DDD removes \emph{minimal underestimating chains}.  A minimal underestimating chain is an underestimating chain, $(\w_1,...,\w_k)$ for which $(\w_1,...,\w_{k-1})$ and $(\w_2,...,\w_{k})$ correspond to resource feasible paths.  Start by considering a VRP where the only disposable resource is time $t$ like the PDPTW.  For each sub-chain $(\w_1,...,\w)$ of the minimal underestimating chain, a timed node is added corresponding to the true end time.  The timed node is $((\f_\w^-, (\tau_{\f_\w}^r)_{r\in R_1}), \tau^t_{\f_\w}\circ...\circ\tau^t_{\f_{\w_1}}(\alpha))$, where $\alpha$ is the earliest start time at the beginning of the represented path.  Timed fragments are added to leave the new timed nodes, and the end timed nodes of existing timed fragments are updated if they are no longer the latest candidate.  This is so that the early arrival and longest arc properties are upheld.  Since $\w_1,...,\w_{k-1}$ will necessarily have there end timed nodes updated, the underestimating chain no longer exists.  
Also, no equivalent underestimating chain exists.  This is because no timed copy of $\f_{\w_\kappa}$ leaves the final added timed node, $\vv$, otherwise the chain would be resource feasible, and the longest arc property ensures their is no equivalent chain to $(\w_1,...,\w_{k-1})$ that arrives at a resourced node earlier than $\vv$.

If there are multiple disposable resources we may need to add more resourced nodes than described above.  For example, suppose the VRP considers load $q$ as well as time $t$ and that end resource consumptions are rounded to the nearest value of $q$ first in the REN. If an underestimating chain is infeasible with respect to load, then the above process will remove all representations of the underlying infeasible path.  However, if an underestimating chain is infeasible with respect to time, then adding resourced nodes that only consider the minimum possible load at the start of the represented path may fail to remove equivalent underestimating chains that start with a higher load.  We therefore calculate $\V(\w_1,...,\w_k)$ the set of resourced nodes in the network where an equivalent underestimating chain can start from.  Then, we must add resourced nodes that correspond to each load $q$ that a resourced node in $\V(\w_1,...,\w_k)$ has.

\subsection{Enhancement 2: Formulation Leveraging}\label{formlev}
Our second proposed enhancement is to use a well known variable fixing technique of BPC algorithms.  It alleviates the relaxation-size trade-off by reducing the size of the fragment formulation through filtering out many resourced fragments that cannot be in an optimal solution.  The technique relies on the following well known theorem that is used in many applications for solving IPs. 
\begin{theorem}\label{varfix}
Suppose we have an IP, $P:=\{\min c^\top x: Ex=b, x\in \mathbb{Z}^d_+\}$ and an upper bound, $z_{ub}$.  The LP $D:=\{\max y^\top b : y^\top E \leq c\}$ is the dual problem of the linear relaxation of $P$.  Given a feasible solution to $D$, $y^*$, the reduced cost of variable $x_j$ is $\bar{c}_j=c_j - y^{*\top} E_j$ for each $j=1,...,d$ where $E_j$ is the column in $E$ corresponding to variable $x_j$.  The objective value of $D$ at $y^*$ is a valid lower bound on $P$, $z_{lb} = y^{*\top}b$.  Given variable $x_j$ for some $j\in\{1,...,d\}$, if $\bar{c}_j > z_{ub} - z_{lb}$ then $x_j = 0$ in all optimal solutions of $P$.
\end{theorem}

This result is used directly on the fragment formulation by \cite{frag3} for the VRP they focus on.  However for the PDPTW and possibly other VRPs, the linear relaxation of RFN does not possess a strong enough lower bound and $z_{ub} - z_{lb}$ is too large for $|\W|$ to be reduced sufficiently.  The route formulation can have a significantly stronger linear relaxation than RFN especially when augmented with valid inequalities.  However, enumerating all routes with reduced cost less than this smaller gap is often still intractable.  The overarching idea of our enhancement therefore, is to use the strong lower bound from the root MP to significantly reduce the number of resourced fragments instead.  We coin the name \emph{formulation leveraging} (FL) for this concept, seeing as we are leveraging the lower bound of a strong formulation to make a weaker formulation solvable.

The process uses the well known variable fixing technique from \cite{vf2}.  Typically, edges $(i,j)$ are removed from $G$ by calculating a lower bound on $\rho_{ij}$, the minimum reduced cost of any route containing the edge.  In BPC, the aim is to increase the efficiency of the labelling algorithm.  On the other hand, we use it to remove resourced fragments $\w$ by calculating a lower bound on $\rho_\w$, the minimum reduced cost of any route with a representation that contains $\w$.  If this value is greater than the optimality gap of SPP then $\w$ is removed because it cannot be in a representation of any route in an optimal solution.

To find a lower bound on $\rho_\w$, MP must be solved by CG.  The technique uses undominated forward and backward labels that arise from solving the final PP with either forward or backward labelling to construct resource feasible (possibly non-elementary) paths from the start depot to the end depot that include the path of resourced fragment $\w$ as a sub-path.  As Theorem \ref{varfix2} states, the reduced cost of at least one of these paths is guaranteed to be a lower bound on $\rho_\w$.  In our description of FL, we assume the PP is solved as an SPPRC rather than an ESPPRC as this is the case in our implementations for the VRPs discussed in this paper.  The modifications to the results in this section are straightforward in the case where the PP is solved as an ESPPRC. 

After using forward labelling we have $\FF(\f)$, the set of undominated forward paths that a vehicle can complete directly before completing fragment $\f\in\F$.  This is the set of undominated forward paths $\p$ that end with non-disposable resource consumptions equal to those of $\f$, $L^1(\p) = L^1_\f$, and disposable resource consumption vector in the feasible set for fragment $\f$, $L^2(\p)\in T_\f$.  In the PDPTW, these are the forward paths with the same set of open deliveries as fragment $\f$, that arrive early enough at the first request of $\f's$ path for the vehicle to then complete the fragment.  After using backward labelling we have $\BB(\p)$, the set of undominated backward paths that can be completed by a vehicle after forward path $\p\in\P$.  This is the set of undominated backward paths $\bp$ with non-disposable resource consumption $\bL^1(\bp)$ that is compatible with $L^1(P)$ and disposable resource consumption vector that satisfies $\bL^2(\bp) \geq L^2(\p)$ component-wise.  In the PDPTW, these are the backward paths that have the same onboard delivery set as forward path $\p$ and leave the head vertex of $\p$ later than the vehicle arrives after traversing the forward path. Note that for the PDPTW, forward and backward non-disposable resource consumptions are only compatible if they are equal.  In general, the definition of compatibility is problem specific.  Finally, let $\bar{c}_{\w}$ be the reduced cost contribution of resourced fragment $\w$ to any route it is in a representation of.  That is, given $\p_{\w} = (i_1,...,i_k)$, $\bar{c}_\w = \sum_{j=1}^{k-1}(c_{i_j,i_{j+1}}-\pi_{i_j}/2 - \pi_{i_{j+1}}/2)$. 
\begin{theorem}\label{varfix2}
    Suppose we have optimal dual values $\pi = (\pi_i)_{i\in V'}$ to the root MP.  Then we have
    \begin{equation}
        \rho_\w' = \min_{\substack{\p\in\FF(\f_\w) \\ \bp\in\BB(\p\oplus \p_{\w})}}\bar{c}_{\p} + \bar{c}_\w + \bar{c}_{\bp}\leq\rho_\w\label{varfixeq}
    \end{equation}
    for each resourced fragment $\w\in\W$.
\end{theorem}

Theorem \ref{varfix2} can be used directly on fragments because it only depends on $\f_\w$.  We use it to filter out fragments before resourced copies are made. One can calculate a better lower bound on $\rho_\w$ for the resourced fragments that do eventuate using their resource consumption information.  The idea is to only consider forward and backward paths that have disposable resource consumption that is congruent with the start and end disposable resource consumption of the resourced fragment.  In order to formally define this, let $\Lambda_\vv$ be the set of true end disposable resource vectors a resourced fragment that ends at $\vv$ or an earlier resourced node can have.  If there is one disposable resource like in the PDPTW, then each $L^2\in\Lambda_\vv$ satisfies $L^2<L^2_{\vv'}$ for all later resourced nodes $\vv'\in\mathtt{later}(\vv)$.  Otherwise, the set depends on the rounding order of the disposable resources.  

\begin{theorem}\label{varfix3}
    Suppose without loss of generality that all representations $(\w_1,...,\w_k)$ must satisfy $L^2(\p_{{\w_1}}\oplus...\oplus \p_{{\w_j}})\in \Lambda_{\w_{j+1}^+}$ for $j=1,...,k-1$ and $\tau^2_{\w_j}\leq L^2(\p_{{\w_1}}\oplus...\oplus \p_{{\w_j}})$ for $j=1,...,k-1$. Any representation can be modified to an equivalent one of these because of the early arrival property.  In this case,
    \begin{equation}\label{varfixeq2}
    \begin{aligned}
        \rho_\w'' &= \min_{\substack{\p\in\FF(\f_\w)\text{ s.t. } L^2(\p)\in \Lambda_{\w^+}\\ \bp\in\BB(\p\oplus \p_{\w})\text{ s.t. }\bL^2(\bp)\geq \tau^2_\w}}\bar{c}_{\p} + \bar{c}_\w + \bar{c}_{\bp} \\
        &\leq\rho_\w,
    \end{aligned}
    \end{equation}
    for each resourced fragment $\w\in\W$.
\end{theorem}

\subsubsection{Rank-1 Cuts}
The route formulation's optimality gap can be reduced by adding valid inequalities to MP which increases the effect of FL.  Rank-1 cuts (R1Cs) are one family of valid inequalities first applied to VRPs by \cite{ssr}.  Given a request set $C\subset V'$ and a multiplier $p_i$ for each $i\in C$ the inequality is of the form 
\begin{equation}
    \sum_{\rr\in\R}\bigg{\lfloor}\sum_{i\in C} a_\rr^ip_i\bigg{\rfloor} \lambda_\rr \leq \bigg{\lfloor} \sum_{i\in C}p_i\bigg{\rfloor}.
\end{equation}
Such cuts are usually strong, but significantly increase the complexity of solving PP because each added cut increments the dimensionality of the labels for dominance purposes.  See \cite{ssr} for details.  \cite{lmssr} introduce the limited memory technique which reduces this impact.  We use limited memory R1Cs (lm-R1Cs) with $|C|$ equal to 1,3,4, and 5.  See \cite{lmssr} for the multiplier values and how lm-R1Cs affect solving PP.  It is straightforward to modify the calculation of $\rho_\w'$ and $\rho_\w''$ given the incorporation of lm-R1Cs.

\subsection{Enhancement 3: Column Enumeration for Row Elimination}
\label{cere}
The effective filtering of FL provides the foundation for a third enhancement that strengthens RFN by lengthening the fragment paths.  As previously mentioned, doing so causes an exponential increase in the number of variables. However, many of the new variables can be filtered out using FL.  Because of the particular way fragments are lengthened, the enhancement also reduces the number of flow conservation constraints (\ref{rfnc3}).  Therefore, the final enhancement, column enumeration for row elimination (CERE), alleviates the relaxation-size trade-off by strengthening RFN and reducing the number of constraints, without risking the number of variables becoming too large.

The idea is to concatenate fragment paths end-to-end.  Given a node $\eta\in \HH$, CERE enumerates all fragments that can be constructed by joining a fragment that ends at $\eta$ and one that starts at $\eta$.  Let $\F_\eta^+$ and $\F_\eta^-$ be the set of fragments that start and end at $\eta$ respectively.  Let $\F_\eta^\cup$ be the set of fragments obtained by joining each in $\F_{\eta}^-$ with each in $\F_{\eta}^+$.  That is, each $\f\in\F_\eta^\cup$ has path $\p_\f = \p_{\f_1}\oplus \p_{\f_2}$ and $L^1_\f = L^1_{\f_1}$ for some $\f_1\in\F_\eta^-$ and $\f_2\in\F_\eta^+$.  Clearly $\F_\eta^\cup$ can replace $\F_\eta^+$ and $\F_\eta^-$ and each route will still have a representation in the network. In the event that the replacement occurs, we say node $\eta$ is absorbed by CERE.  This is because no fragment in $(\F\cup\F_\eta^\cup)\setminus(\F_\eta^-\cup\F_\eta^+)$ starts or ends at $\eta$.  Thus, no resourced copies of $\eta$ are needed and $|\V|$ is reduced.

We only perform CERE on nodes $\eta$ that do not cause a significant increase in the number of fragments.  For each node $\eta\in\HH$ we enumerate set $\F_\eta^\cup$ and filter out members using FL.  The node is absorbed if the number of fragments increases less than a chosen integer $\mathtt{MaxIncrease}$.  This parameter controls the aggressiveness of CERE.  Larger values result in more fragments with longer paths and setting it sufficiently large will result in route enumeration and the fragment formulation simplifies to its route counterpart.  Choosing moderate values strengthens RFN and reduces the number of constraints significantly, without the number of variables becoming too large.

\subsection{The Overarching Strategy}
While the three enhancements are intended for improving the fragment approach on the PDPTW, they lead to a general framework that may theoretically work on VRPs other than PDPs. For instance, the combination of FL and CERE fixes the problem of fragment paths being too limited in length for RFN to be sufficiently strong.  Constructing the network with DDD enables the fragment approach to handle problems with wider resource windows at each request, since resourced nodes are added sparingly.  In essence, the approach offers a solution method that is an alternative to BPC, which is theoretically applicable to many VRPs.  Whilst BPC seeks to solve a strong formulation even though a bespoke BB framework is invoked and each BB node is computationally expensive to process, the enhanced fragment approach seeks to solve a strong formulation that is additionally constrained to be small enough for a commercial BB solver to handle.  Each of the enhancements contributes to the automatic construction of such a strong formulation.   The advantage of using a commercial BB solver is that more effective branching and cutting plane strategies can be used and significantly larger BB trees can be explored.  Further research is needed to determine whether other VRPs can be more efficiently solved by this new strategy.  We begin this process by applying our approach to the TDDRP additionally to the PDPTW.  We choose the TDDRP because while it is not a PDP, it is another example where route complexities lead to a natural fragment definition.

\subsection{The Overall Approach}
We conclude this section with the steps of the proposed approach.  Both VRPs considered minimize the number of vehicles used and then the travel cost.  We handle these separately.  We target a guess for the minimum number of vehicles $z_{ub}^v$ and increase it if the VRP is infeasible with this target.  Two solve phases are executed for each target number of vehicles.  We guess an upper bound for the travel cost $z_{ub}^c$.  If the optimal travel cost found by solving RFN exceeds this guess, we increase the guess to the cost found and resolve starting from Step \ref{fragenum}.  We choose a small and a large $\mathtt{MaxIncrease}$ value to use for CERE in the first and second phases, $\mathtt{MaxIncrease'}$ and  $\mathtt{MaxIncrease''}$ respectively.  In the second phase, we only attempt CERE on nodes that were removed in the first phase.  Removing all of these is not guaranteed because the travel cost optimality gap is larger.  If removal is unsuccessful for some nodes in the second phase, the first phase must be resolved over the new node set.  This is because the branching constraint added in the second phase in Step \ref{formulate} is only valid if both phases are solved with the same set of nodes.
\begin{enumerate}
    \item Solve MP by CG with the objective set to minimizing the number of vehicles.  Let the optimal objective be $z_{lb}^v$ and a minimum vehicle number guess be $z_{ub}^v = \lceil z_{lb}^v\rceil$.
    \item Add the constraint $\sum_{\rr\in\R'}\lambda_\rr = z_{ub}^v$ to MP and solve it by CG with the objective set to minimizing the travel cost.  Let $z_{lb}^c$ be the optimal objective value and $z_{ub}^c > z_{lb}^c$ be an upper bound guess for the optimal travel cost.  Let $z_{force}^c = -\infty$.\label{tclb}
    \item Enumerate fragments $\F$ and perform FL according to Theorem \ref{varfix2} to reduce $|\F|$.  Both optimality gaps $z_{ub}^v - z_{lb}^v$ and $z_{ub}^c-z_{lb}^c$ can be used for this.\label{fragenum}
    \item Perform CERE.  If $z_{force}^c=-\infty$ use $\mathtt{MaxIncrease'}$ and let $\HH'$ be the set of remaining nodes.  Otherwise use $\mathtt{MaxIncrease''}$ and only attempt CERE on nodes not in $\HH'$.  If $\HH''$ is different from $\HH'$, reset $z_{ub}^c$ to its initial value, $z_{force}^c$ to $-\infty$ and $\HH'$\ to $\HH''$.\label{cerestep}
    \item Build $\N$ and use FL according to Theorem \ref{varfix3} to reduce $|\W|$.\label{build}
    \item Formulate RFN.  The objective minimizes the travel cost and the number of vehicles is set to $z_{ub}^v$. If $z_{force}^c > -\infty$, include the branching constraint $\sum_{\w\in\W|\rho_\w''\geq z_{force}^c-z_{lb}^c}x_\w\geq 1$.\label{formulate}
    \item Solve RFN via DDD yielding travel cost $z_{ub}^{RFN}$ which may be $\infty$.  If $z_{ub}^{RFN}\leq z_{ub}^c$ and $z_{ub}^{RFN}<\infty$ then stop with an optimal solution.  Otherwise if $z_{force}^c = -\infty$, set $z_{force}^c = z_{ub}^c$, $z_{ub}^{c}= z_{ub}^{RFN}$ and go to Step \ref{fragenum}.  Otherwise increment the target number of vehicles guess $z_{ub}^v$ and go to Step \ref{tclb}.
\end{enumerate}

The resourced nodes added by DDD are reused each time the resource expanded network is built in Step \ref{build}.  We discuss the choice of initial resourced nodes in the next sections.   

\section{Computational Results for the PDPTW}\label{pdptw_old}
Two instance sets have been used when testing exact PDPTW algorithms: those from \cite{pdptw3} and \cite{pdptw4}.  Those commonly considered from \cite{pdptw3} were all solved by \cite{pdptw2} using route enumeration and variable fixing after solving MP by CG.  This is a special case of our approach where CERE is performed on all nodes.  Some \cite{pdptw4} instances cannot be solved by this method and either BPC or our fragment approach with selective CERE is necessary.  We consequently focused on these.  They are separated into classes AA, BB, CC and DD.  Classes AA and BB have narrower time windows than CC and DD.  Classes AA and CC have smaller vehicle capacity than BB and DD.

We implemented our algorithm in Rust and solved RFN using Gurobi 11.  The computer had an Intel i7-12700H 2.70 GHz processor with 64GB of memory.  \cite{bapcod}'s BaPCod in conjunction with \cite{vrpsolver}'s VRPSolver were locally installed and used to solve MP by CG with CPLEX Studio 22.1.1. The configuration we chose for completing steps 1 and 2 of the overall algorithm is outlined in the online supplement.  The configuration when solving purely with VRPSolver was unchanged from \cite{vrpsolver}.

We run six variants of our algorithm that each include different combinations of the enhancements.  Table \ref{variants} summarises the variants.  Observe that CERE is never used without FL as doing so is impractical.  When DDD is not used, RFN is solved by BC instead.

\begin{table}[!htp]\centering
\caption{Variants of the fragment algorithm.}\label{variants}
\scriptsize
\begin{tabular}{lrrrr}\toprule
Variant &DDD &FL &CERE \\\midrule
B &- &- &- \\
D &\checkmark &- &- \\
BF &- &\checkmark &- \\
DF &\checkmark &\checkmark &- \\
BFC &- &\checkmark &\checkmark \\
DFC &\checkmark &\checkmark &\checkmark \\
\bottomrule
\end{tabular}
\end{table}

For each instance, we use the heuristic upper bounds from \cite{pdptw4}, $z^{heuristic}$, as initial best known solution costs and choose $z_{ub}^c = \min\{z^\text{heuristic}, z_{lb}^c+\psi\}$ where $\psi$ is 20 if the best known heuristic solution has the same number of vehicles as targeted and 30 if not.  When solving with BC, $\N$ has resourced copies of each node that are 5 units of time apart.  When solving with DDD, $\N$ starts with a single resourced copy of each node at the earliest time possible.

\begin{figure}[h]

{\includegraphics[width=0.5\textwidth]{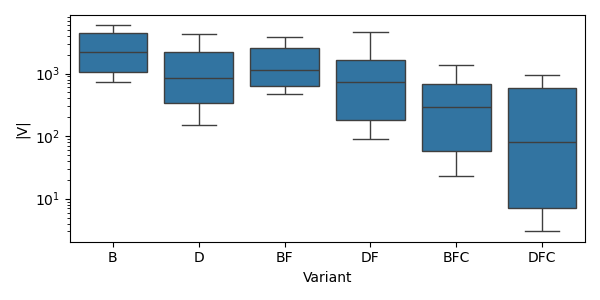}}
\hfill
{\includegraphics[width=0.5\textwidth]{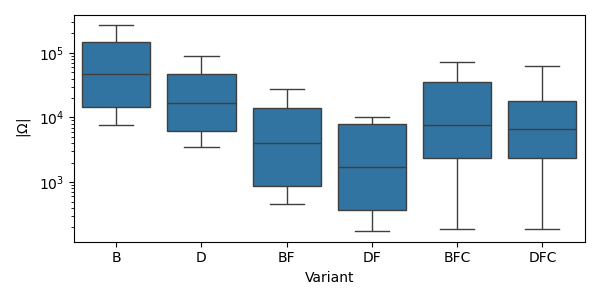}}

\caption{Number of timed nodes and fragments in the final network across the PDPTW instances.}
\label{pdptw_size}
\end{figure}

Figure \ref{pdptw_size} reports boxplots of the number of timed nodes and timed fragments in the final resource expanded network across the benchmark instances that cannot be solved at the root node by CG. Those that are, solve in 21 seconds on average. Comparing the variants that use DDD to those that use BC shows that DDD enables fewer timed nodes in the network, and consequently fewer timed fragments.  Higher variance is observed in the DDD variants because of the adaptive nature of DDD.  Some instances require a lot of timed nodes to be added whilst other require very few.  As expected, using FL without CERE significantly reduces the number of timed fragments in the network.  Finally, CERE reduces the number of timed nodes and increases the number timed fragments without the number becoming too large.  

\begin{figure}[h]
{\includegraphics[width=\textwidth]{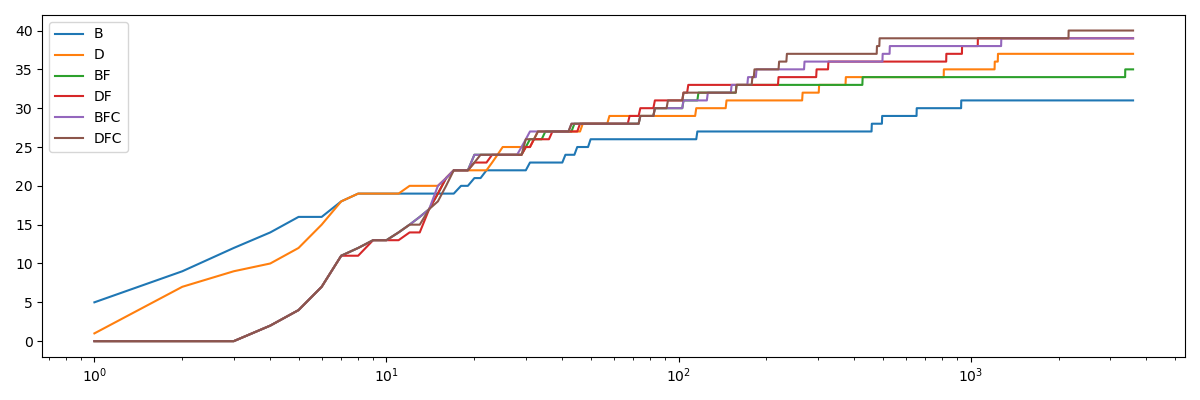}}
\caption{Number of instances solved to optimality within a given amount of time.
}
\label{pdptw_solve_times}
\end{figure}

Figure \ref{pdptw_solve_times} compares the performance of the variants across all benchmark instances by reporting the number of instances solved to optimality within a given amount of time.  The variants that implement FL outperform those that do not. Variants B and D solve the easier instances more efficiently because of the extra time needed by FL to solve MP by CG and perform variable fixing.  However, the reduced number of timed fragment variables resulting from FL enables the harder instances to be solved. Each DDD variant outperforms its BC counterpart.  In fact, variant D even outperforms BF, despite the latter generally solving a formulation with fewer timed fragment variables. This is because for some instances, the static network of BF does not include enough resourced copies of some nodes, leading to the addition of too many feasibility cuts which slows convergence.  Variant D on the other hand adds timed nodes as needed and does not have to use feasibility cuts.  Observe that BFC outperforms BF because lengthening the fragments with CERE means fewer underestimating chains exist in the network, leading to fewer feasility cuts.  The underestimating chains that are removed are generally shorter, because fragment paths are longer, and the feasibility cuts are more effective for such chains. As expected, the best performing variant includes all enhancements. 

Finally, Table \ref{compvrpsolver} shows that the fragment method is competitive with, and can outperform \cite{vrpsolver}'s BPC algorithm.  We focus on the eight hardest benchmark instances.  The other 32 are solved by the fragment method in 26 seconds and by BPC in 14 seconds on average - times that are negligible compared to those reported in Table \ref{compvrpsolver}.  The total number of nodes explored in DDD iterations is reported for the fragment algorithm, using either one or eight threads in the commercial solver.  As shown, the fragment method outperforms BPC on the majority of these instances when using one thread.  The only instance where BPC significantly outperforms the one-threaded fragment method is DD60.  Here, the travel cost lower bound of the MP is $z^{c}_{lb} = 2010.89$ and the optimal solution cost is 2069.175.  The fragment method easily finds this optimal solution, however, the MP optimality gap is relatively large compared to other instances, and FL and CERE are less effective at constructing a strong and small fragment formulation.  Consequently, many BB nodes must be explored in the final DDD iterations to prove optimality.  In BPC, the MP lower bound is improved significantly to 2031.31 after branching once.  Thus, it is likely that the fragment method would be efficient with the smaller MP optimality gaps of the resulting child nodes.  Such augmentation of BPC and the fragment method may be an effective acceleration technique to explore in future research.  Even without this extension, if we run Gurobi with 8 threads as in \cite{frag1}, we solve DD60 in under one hour.  Having instant access to all of Gurobi's acceleration techniques including parallelization is part of the motivation for the fragment method.  Parallelizing BPC on the other hand requires extra effort, with an unclear benefit to solve time.

\begin{table}[h!]
\centering
\footnotesize
\caption{Comparison with \cite{vrpsolver}'s BPC.}
\label{compvrpsolver}
\begin{tabular}{l|rr|rrrr}\toprule
&\multicolumn{2}{c|}{Nodes} &\multicolumn{3}{c}{Solve Time (s)} \\\cmidrule{2-6}
Instance &1 Thread &8 Threads &1 Thread &8 Threads &BPC \\\midrule
CC60 &6227 &5457 &170.5 &144.7 &362.8 \\
CC65 &8903 &7301 &319.7 &187.4 &562.6 \\
CC70 &5373 &13655 &555.4 &585.9 &383.4 \\
CC75 &613 &541 &105.5 &98.4 &186.8 \\
DD60 &35239 &80468 &3600.0 &2530.2 &1435.3 \\
DD65 &10220 &5360 &345.1 &178.4 &786.4 \\
DD70 &4647 &5629 &330.4 &261.7 &689.6 \\
DD75 &4845 &4522 &659.8 &534.0 &3016.1 \\
\bottomrule
\end{tabular}
\end{table}

\section{The TDDRP}\label{tddrp}
In this section we extend the framework to the TDDRP proposed by \cite{tddrp}.  We deviate from the route, fragment and fragment REF definitions in Section \ref{algorithm} because TDDRP routes are not only truck paths, but consist of a truck and a drone working in tandem to deliver requests.  Each drone can be carried by its truck or deployed to service some requests and recollected at another location.  We have graph $G=(V,E)$ with $V = \{0,n+1\}\cup N$ where $N = \{1,...,n\}$ is the set of requests.  Each request $i$ has a time window $[a_i, b_i]$ and a weight $d_i$.  We assume an infinite fleet of homogeneous truck-drone pairs.  The drone of each pair has limited flying duration $\mathtt{L}^\text{drone}$ caused by its battery capacity.  The drone must return to the truck to exchange its battery before reaching its maximum flying duration.  We let $Q^\text{drone}$ and $Q^\text{truck}$ be the drone and truck capacities respectively with $Q^\text{drone}<Q^\text{truck}$.  Each edge $(i,j)\in E$ has a mileage cost and travel time for the truck and drone, $c^\text{truck}_{ij}$ and $t^\text{truck}_{ij}$, and $c^\text{drone}_{ij}$ and $t^\text{drone}_{ij}$ respectively.  Note that $c^\text{drone}_{ij}<c^\text{truck}_{ij}$ and $t^\text{drone}_{ij} < t^\text{truck}_{ij}$ for all $(i,j)\in E$.  Each truck-drone pair has sufficiently large fixed cost $F$ such that the TDDRP aims to minimize the number of routes and then the mileage cost of covering each request exactly once.

\subsection{Drone Flights}
Our approach is based on enumerating all paths that are feasible for a drone to complete.
\begin{definition}
A \emph{drone flight} is an elementary path $\f^\text{drone} = (i_1,...,i_\kappa)$ in $G$ that respects the drone's maximum flight time and capacity assuming that it is deployed at $i_1$ and collected at $i_k$.  

It is possible to enumerate all drone flights because their battery and capacity are small.  Let $d_{\f^\text{drone}}=\sum_{j=2}^{k-1}d_{i_j}$ be the total weight of requests delivered by the drone. The drone must reach requests $i_2$ to $i_{k-1}$ within their time windows.  Let $s_i$ be the earliest arrival at each request $i=i_1,...,i_k$.  So $s_{i_1} = t_{0i_1}^\text{truck}$, $s_{i_2} = s_{i_1} + t_{i_1i_2}^\text{drone}$ and $s_{i_{j+1}} = \max\{a_{i_j}, s_{i_j}\} + t_{i_1i_2}^\text{drone}$ for $j=2,...,k-1$.  We require $s_{i_j}\leq b_{i_j}$ for $j=2,...,k-1$.  Let $\mathtt{latest}_{\f^\text{drone}}$ be the latest deployment time for the drone to complete the flight feasibly.  For each deployment time $t^*\in [t_{0i_1}^\text{truck}, \mathtt{latest}_{\f^\text{drone}}]$, let $\mathtt{finish}_{\f^\text{drone}}(t^*)$ be earliest arrival at $i_k$.  Let $\F^\text{drone}$ be the set of drone flights and $\F^\text{drone}_i$ be the drone flights deployed at $i\in V$.  
\end{definition}

The labelling algorithm for calculating negative reduced cost routes is similar to \cite{tddrp}'s.  Three differences stem from using drone flights.  Firstly, our labels have a different form, $\l = (i_\text{truck}, i_\text{drone}, \bc, (t_\text{truck}, t_\text{drone}, q))$, where $i_\text{truck}\in V$ and $i_\text{drone}\in V$ are the truck and drone locations respectively, $t_\text{truck}$ and $t_\text{drone}$ are the truck and drone times respectively, and $q$ is truck capacity used by requests delivered.  We do not track the drone capacity used like \cite{tddrp} since this is satisfied by every drone flight.  Secondly, when $i_\text{truck}=i_\text{drone}$, only two extensions are possible. In forward labelling, either the truck carries the drone to a new location, or a drone flight is deployed.  In backward labelling, either the truck carries the drone or leaves it at a location to be deployed to a future truck location. Finally, when $i_\text{truck}\neq i_\text{drone}$, the label is only extended by moving the truck in forward labelling.  In backward labelling, either the truck moves again, or the drone completes a drone flight to rendezvous with the truck.  Full details are given in Appendix 4.

\subsection{Routes and Fragments}
Routes are complex to directly define for the TDDRP.  Instead, we define them as a chain of fragments.  Each fragment consists of a path for the drone and the truck.  When a drone flight is deployed, the vehicle completes a \emph{truck leg} which finishes at the same location as the drone flight.
\begin{definition}
Given a drone flight $\f^\text{drone} = (i^1_1,...,i^1_{\kappa_1})$, a \emph{truck leg} is an elementary path $\f^\text{truck} = (i^2_1,...,i^2_{\kappa_2})$ such that the total weight of requests serviced by $\f^\text{drone}$ and $\f^\text{truck}$ does not exceed $Q^\text{truck}$, each location in $\f^\text{truck}$ can be reached by the truck within its time window, and $\f^\text{truck}$ only shares its start and end locations in common with $\f^\text{drone}$.  

Let $\mathtt{latest}_{\f^\text{truck}}$ be the latest start time at $i_1^2$ for which each request can be reached in its time window. For every start time $t^*\in[t_{0i_1^2}^\text{truck}, \mathtt{latest_{\f^\text{truck}}}]$, let $\mathtt{finish}_{\f^\text{truck}}(t^*)$ be the earliest arrival at $i^2_{k_2}$. Finally, let $d_{\f^\text{truck}}=\sum_{j=2}^{k_2}d_{i_j^2}$ be the total weight of requests delivered by the truck leg.
\end{definition}

\begin{definition}
A \emph{TDDRP fragment} is a drone path $\f^\text{drone}$ with a truck path $\f^\text{truck}$ where either:
\begin{enumerate}
    \item The pair represents the truck carrying the drone, $\f^\text{drone} = \f^\text{truck} = (i,j)\in E$.
    \item Drone path $\f^\text{drone}$ is a drone flight and $\f^\text{truck}$ is a corresponding truck leg. 
\end{enumerate}  
We refer to these as type 1 and type 2 fragments respectively.
\end{definition}

Enumerating type 2 fragments can be intractable because truck legs are loosely resource constrained.  We circumvent this by using FL according to Theorem \ref{varfix2} on partial fragments consisting of a drone flight and a partial truck leg.  Details are in the online supplement.

\subsection{Resource Expanded Network for the TDDRP}

A truck-drone pair synchronizes at the start and end of every fragment.  So nodes need only one location and the only resources the network must consider are the cumulative load of the truck-drone pair $q$ and the earliest time that the pair have both arrived $t$. So for each type 1 fragment $\f = ((i,j),(i,j))$, define fragment REFs $\tau_\f^t(t^*) = \max\{a_i, t^*\} + t_{ij}^\text{truck}$ and $\tau_\f^q(q^*) = q^* + d_j$. For each type 2 fragment $\f = (\f^\text{drone}, \f^\text{truck})$, define $\tau_\f^q(q^*) = q^* + d_{\f^\text{drone}} + d_{\f^\text{truck}}$ and $\tau_\f^t(t^*) = \max\{\mathtt{finish}_{\f^\text{drone}}(t^*), \mathtt{finish}_{\f^\text{truck}}(t^*)\}$.  Resources $t$ and $q$ are disposable.  So $R^1 = \emptyset$ and $R^2 = \{t, q\}$.  Type 1 fragments $\f=((i,j),(i,j))$ have $T_\f = [t_{0i}^\text{truck}, b_j-t_{ij}^\text{truck}\}]\times [0,Q^\text{truck}-d_j]$ and type 2 fragments $\f = (\f^\text{drone}, \f^\text{truck})$ have $T_\f = [t_{0i}^\text{truck}, \min\{\mathtt{latest}_{\f^\text{drone}},\mathtt{latest}_{\f^\text{truck}}\}]\times [0,Q^\text{truck}-d_{\f^\text{drone}}-d_{\f^\text{truck}}]$.  The form of resourced nodes and fragments then follows from Section \ref{ren}.

\subsection{Computational Results}
The instances solved in \cite{tddrp} are not publicly available.  So we test the performance of our method on instances generated according to their description.  These have $Q^\text{truck}=100$ and $Q^\text{drone}=20$.  Routes must finish before time $480$.  The average request weight is $30$ and the average time window width is $60$. Travel times satisfy $t_{ij}^\text{truck} = 2t_{ij}^\text{drone}$ and the truck travel cost per minute is approximately three times that of drones.  Cheaper drones than trucks encourages using fragments with short truck legs, making FL important for filtering out those with long truck legs.

Instances with clustered requests and a centred depot, with uniformly random request and depot coordinates, and with uniformly random request coordinates and a centred depot are generated.  \cite{tddrp} alter the proportion of requests that can be serviced by the drone between 10\% and 90\%.  We only consider the hardest case of 90\%.  Three instances per type are generated for each $n\in\{20,35,45,65,85,105,125\}$.  The instances can be found at github.com/ls1704/tddrpdata.

We use our own CG implementation to solve MP.  Three rounds of cut separation are used in step 2 of the overall algorithm, each time resolving MP to optimality.  No cuts are separated in step 1 as preliminary testing showed that a solution with $\lceil z_{lb}^v\rceil$ vehicles could be found for every instance.  We set $z_{ub}^c = \min\{z^\text{RMP}, z_{lb}^c+0.5\}$ where $z^\text{RMP}$ is the cost of the optimal integer solution of RMP after completing steps 1 and 2.  When performing CERE, we set $\mathtt{MaxIncrease'}=500$ and $\mathtt{MaxIncrease''}=\infty$.  We initialize $\N$ with resourced copies of each request $i\in V'$ of the form $(i, 10\delta, t^\text{truck}_{0i})$ for $\delta\in\{d_i/10,(d_i+10)/10,...,(Q^\text{truck}-10)/10,Q^\text{truck}/10\}$.

Table \ref{tddrpres} summarises the performance of our method.  Each row gives average results for the nine corresponding instances.  Column Gap gives the initial optimality gap between $z_{lb}^c$ and the optimal cost.  Column Nodes gives the proportion of nodes removed from $\HH$ by CERE.  Column Optimal gives the number of instances solved to optimality within the time limit.  Columns CG, $|\W|$, $|\V|$, and Solve (s) give the time for MP to be solved by CG in seconds, the final number of resourced fragments, the final number of resourced nodes, and the total solve time in seconds respectively.  The averages are only calculated over those instances solved within the time limit.  We compare the solve time and number of instances solved with the results in \cite{tddrp} via sub-columns Proposed and BPC respectively.

Firstly, Table \ref{tddrpres} shows that the optimality gap $(z_{ub}^c-z_{lb}^c)/z_{lb}^c$ remains small for all instance sizes.  \cite{tddrp} report larger gaps for their instances, however, this is because they solve with a composite objective which multiplies the number of routes by a large fixed cost and then adds the mileage cost.  The number of routes can be fractional, causing the gap to be larger with respect to this objective.  Columns Nodes and $|\W|$ show that CERE is effective at removing the majority of nodes from $\HH$ while FL is effective at keeping the number of resourced fragments tractable.  Both depend on the tight initial optimality gap. Despite there being two disposable resources, the number of resourced nodes remains small for all $n$.  Analyzing columns CG and Solve shows that CG is the limit of our approach rather than solving RFN.  Making the separation of lm-R1Cs less aggressive is unlikely to improve performance, as this would result in larger gaps.  Our approach only fails to solve one instance with 125 requests within one hour.  In comparison, the BPC algorithm cannot solve all instances with 35 and 45 instances.  The average solve times of our approach are also significantly smaller than those in \cite{tddrp}.  Even though solving different instances generated in the same way is not ideal for comparison of performance, these results are strong enough to conclude that our algorithm outperforms the BPC algorithm.

\begin{table}
\centering
\caption{Computational results for TDDRP instances.}
\label{tddrpres}
\footnotesize
\begin{tabular}{lrrrrrrrrr}
\cline{7-10}
 & \multicolumn{1}{l}{} & \multicolumn{1}{l}{} & \multicolumn{1}{l}{} & \multicolumn{1}{l}{} & \multicolumn{1}{l}{} & \multicolumn{2}{c}{Solve (s)} & \multicolumn{2}{c}{Optimal} \\ \hline
\multicolumn{1}{l}{$n$} & \multicolumn{1}{r}{Gap} & \multicolumn{1}{r}{Nodes} & \multicolumn{1}{r}{$|\W|$} & \multicolumn{1}{r}{$|\V|$} & CG (s) & \multicolumn{1}{r}{Proposed} & \multicolumn{1}{r}{BPC} & \multicolumn{1}{r}{Proposed} & \multicolumn{1}{r}{BPC} \\ \hline
\multicolumn{1}{l}{20} & \multicolumn{1}{r}{0.0041} & \multicolumn{1}{r}{1.00} & \multicolumn{1}{r}{93} & \multicolumn{1}{r}{0} & 2 & 2 & 17 & 9 & 9 \\
\multicolumn{1}{l}{35} & \multicolumn{1}{r}{0.0060} & \multicolumn{1}{r}{0.99} & \multicolumn{1}{r}{870} & \multicolumn{1}{r}{2} & 12 & 13 & 2694 & 9 & 7 \\
\multicolumn{1}{l}{45} & \multicolumn{1}{r}{0.0040} & \multicolumn{1}{r}{0.98} & \multicolumn{1}{r}{1901} & \multicolumn{1}{r}{9} & 22 & 24 & 7067 & 9 & 6 \\
\multicolumn{1}{l}{65} & \multicolumn{1}{r}{0.0032} & \multicolumn{1}{r}{0.91} & \multicolumn{1}{r}{10130} & \multicolumn{1}{r}{47} & 128 & 136 & - & 9 & - \\
\multicolumn{1}{l}{85} & \multicolumn{1}{r}{0.0025} & \multicolumn{1}{r}{0.88} & \multicolumn{1}{r}{17860} & \multicolumn{1}{r}{86} & 306 & 391 & - & 9 & - \\
\multicolumn{1}{l}{105} & \multicolumn{1}{r}{0.0036} & \multicolumn{1}{r}{0.77} & \multicolumn{1}{r}{98560} & \multicolumn{1}{r}{191} & 785 & 884 & - & 9 & - \\
\multicolumn{1}{l}{125} & \multicolumn{1}{r}{0.0019} & \multicolumn{1}{r}{0.78} & \multicolumn{1}{r}{86526} & \multicolumn{1}{r}{234} & 1901 & 2112 & - & 8 & - \\ \hline
\end{tabular}
\end{table}

Admittedly, \cite{tddrp}'s BPC algorithm would likely perform better if it employed route enumeration, which may be possible in many instances due to the tight optimality gap and relatively few requests per route.  Note that our implementation would likely still outperform theirs due to our faster CG which uses drone paths and lm-R1Cs instead of R1Cs.  However, it is likely that the fragment approach does not provide significant advantage over route enumeration for these instances, apart from avoiding memory limit violation in a few with the largest optimality gaps.  We therefore analyse the performance of our approach on 27 new instances that are harder.  We take the initial instances with 45 requests and reduce all request weights to 10. We also increase $Q^{\text{truck}}$ to 120 and 140.  Both of these alterations increase the average number of requests per route and consequently, the time to complete CG.  Preliminary testing revealed that using lm-R1Cs with three or more requests made pricing too inefficient, and so we only use those with one request.  This results in looser optimality gaps and thus we set $\mathtt{MaxIncrease'} = 2500$ in order to increase the frequency of absorbing nodes.  We also increase the upper bound guess $z_{ub}^c$ to $\min\{z^\text{RMP}, z_{lb}^c+1.0\}$ for the clustered instances since these tended to have the largest optimality gaps, but the fewest fragments.  We solve each instance three times and increment the number of cut separation rounds when solving the MP each time.  We also increase the time limit to 15 hours.

Table \ref{tddrpres2} exhibits results for the fragment approach on the harder instances.  The first column gives the instance solved.  The second column gives the truck capacity.  The third column reports the number of cut separation rounds (from zero to two) used in the fastest solve.  Columns IGap and FGap give the initial optimality gap guess, and the final optimality gap with respect to the route formulation.  The rest of the columns are the same as in Table \ref{tddrpres}.  Firstly, notice that additional cut separation rounds significantly increase the CG time, meaning the fastest solve usually uses the fewest cut separation rounds necessary for FL and CERE to be effective via a sufficiently small optimality gap.  Secondly, larger optimality gaps are achieved for these instances, leading CERE to absorb fewer nodes.  This increases the number of resourced nodes in the final REN, however, the number of resourced fragments remains manageable.  Despite this, RFN is still solved efficiently by DDD in many cases.  In fact, for many problems, RFN is solved despite the optimality gap being too large for route enumeration to be possible.  This shows the advantage of the fragment approach over one based solely on route enumeration.  Algorithms which branch before attempting route enumeration would also likely be outperformed by the fragment approach here because the MP is too expensive to resolve at each BB node.  Finally, there exists a few instances where DDD takes a long time to solve RFN relative to the CG time.  In these cases, the truck capacity and optimality gap are largest, meaning CERE is less effective at absorbing nodes to strengthen RFN.  This causes it to be slower to solve in each DDD iteration.

\begin{table}[h!]\centering
\caption{Computational results for the hard TDDRP instances.}
\label{tddrpres2}
\footnotesize
\begin{tabular}{lrrrrrrrrrrr}\toprule
Instance &$Q^{truck}$ &Cut Sep. &CG &IGap &FGap &Node &$|\F|$ &$|\W|$ &$|\V|$ &Solve \\\midrule
45-centre-cluster-1 &100 &2 &1541 &0.022 &0.024 &0.53 &191321 &444212 &585 &3031 \\
&120 &2 &5104 &0.022 &0.019 &0.47 &122348 &295879 &665 &11007 \\
&140 &2 &7638 &0.022 &0.020 &0.40 &152229 &439467 &980 &12282 \\
\midrule
45-centre-cluster-2 &100 &1 &532 &0.023 &0.007 &0.43 &51478 &93029 &478 &852 \\
&120 &1 &1648 &0.022 &0.002 &0.49 &40014 &90634 &716 &2360 \\
&140 &1 &3517 &0.023 &0.004 &0.49 &34137 &98854 &848 &5132 \\
\midrule
45-centre-cluster-3 &100 &1 &989 &0.023 &0.014 &0.56 &54887 &123963 &302 &1244 \\
&120 &2 &7157 &0.002 &0.002 &1.00 &2227 &954 &2 &7192 \\
&140 &1 &8491 &0.008 &0.004 &0.84 &5076 &1979 &179 &8863 \\
\midrule
45-centre-random-1 &100 &1 &678 &0.028 &0.037 &0.47 &23519 &52909 &649 &8457 \\
&120 &0 &1195 &0.030 &0.032 &0.64 &14323 &25402 &523 &1302 \\
&140 &0 &1972 &0.032 &0.047 &0.73 &50309 &256113 &655 &3264 \\
\midrule
45-centre-random-2 &100 &0 &1143 &0.034 &0.030 &0.42 &29779 &36542 &370 &1568 \\
&120 &0 &1744 &0.035 &0.041 &0.44 &50376 &99044 &951 &2619 \\
&140 &1 &14232 &0.036 &0.043 &0.33 &42077 &101855 &1091 &43167 \\
\midrule
45-centre-random-3 &100 &0 &761 &0.031 &0.028 &0.49 &24341 &42773 &366 &918 \\
&120 &0 &1599 &0.032 &0.053 &0.51 &91552 &293640 &999 &3780 \\
&140 &1 &11954 &0.034 &0.050 &0.55 &49504 &181736 &1203 &48275 \\
\midrule
45-random-random-1 &100 &1 &2057 &0.031 &0.007 &0.63 &27492 &44893 &409 &2311 \\
&120 &1 &7236 &0.032 &0.018 &0.60 &28912 &14501 &592 &9239 \\
&140 &1 &18749 &0.033 &0.017 &0.62 &25556 &9855 &518 &27020 \\
\midrule
45-random-random-2 &100 &0 &1120 &0.025 &0.012 &0.34 &33302 &61550 &460 &1365 \\
&120 &0 &2227 &0.027 &0.045 &0.42 &112027 &267316 &1062 &9858 \\
&140 &1 &12638 &0.028 &0.032 &0.46 &46619 &111392 &1385 &38867 \\
\midrule
45-random-random-3 &100 &0 &712 &0.023 &0.033 &0.69 &32878 &68500 &433 &925 \\
&120 &1 &4850 &0.024 &0.034 &0.63 &41313 &89382 &930 &11785 \\
&140 &1 &13076 &0.026 &0.039 &0.66 &50554 &158357 &1067 &44487 \\
\bottomrule
\end{tabular}
\end{table}

These results also provide insight on the kind of instances where the fragment approach is most likely to perform well for other VRPs.  Clearly, instances must permit a sufficiently small route formulation optimality gap such that FL and CERE keep the size of the REN manageable whilst producing a strong enough RFN.  Even then, the fragment approach may not provide an advantage over BPC if the MP is not expensive to solve, and routes can be enumerated easily, as this leads to small BB trees whose nodes can be processed quickly in BPC.  Hence, the fragment approach is most likely to perform well relative to BPC when MP is expensive to solve, and route enumeration is difficult even for very small optimality gaps.  This usually occurs when the average length of routes is large, as in the hard instances for the TDDRP, and those studied for the PDPTW.

\section{Conclusion}\label{conclusion}
We have proposed a new fragment method for solving the PDPTW that is enhanced with DDD, FL and CERE.  Additionally to applying the method to the PDPTW, we present the algorithm generally, which communicates its theoretical applicability to other VRPs.  We support this idea by applying it to the TDDRP where routes have structure that leads to a natural fragment definition.  We show via computational results that our approach outperforms previously proposed fragment methodology, and can also outperform the best performing BPC algorithms for the PDPTW and the TDDRP.

The proposed framework has limitations that inform future research directions.  Firstly, it is difficult to model complex objective function costs in the fragment formulation, like those that depend on combinations of resource fragments or the resource consumption values.  Future research should determine whether the fragment approach can efficiently solve problems that incorporate such costs.  Secondly, for some VRPs, the proposed method may still fail to produce a strong enough fragment formulation, leading to inefficiency in DDD.  Future research should therefore seek other techniques to further strengthen the fragment formulation.  Finally, while we gave some insight into the types of VRP instances that may be solved most efficiently by our method in the last section, it is still unclear whether our method will be effective on other VRPs.  Towards this end, we plan to test its performance on the VRPTW.  \cite{vrpsolver}'s BPC algorithm has all modern enhancements for solving the VRPTW.  So, if our algorithm is competitive with this algorithm for such a thoroughly researched problem, it would strengthen the case for further research into its use.

\section*{Acknowledgments}
Lucas Sippel is supported by an Australian Government Research Training Program Scholarship.

\bibliographystyle{main} 
\bibliography{main}

\end{document}


\maketitle

\section{Proofs} \label{proofs}
\subsection{Theorem 1}
Firstly, the fragment rules ensure an appropriate sequence of fragments $\f_1,...,\f_\kappa$ exists for $\rr$. Properties 1 and 3 imply the existence of resourced fragment $\w_1$ with $\f_{\w_1} = \f_1$ and $\L^2_{\w_1} = (\alpha_{i_1}^r)_{r\in R_2}$.  This combined with Property 3 implies that $\L'^2_{\w_1} \leq \L^2((i_1,...,\f_1^-))$.

For some $j=1,...,\kappa-1$, suppose we have $\w_j$ with $\f_{\w_j}=\f_j$ such that $\L'^2_{\w_j} \leq \L^2((i_1,...,\f_j^-))$.  Since all resources $r\in R_2$ are disposable, $\L'^2_{\w_j}\in T_{\f_{j+1}}$.  Therefore by Property 2, there exists resourced fragment $\w_{j+1}$ with $\f_{\w_{j+1}} = \f_{j+1}$ and $\L^2_{\w_{j+1}} = \L'^2_{\w_j}$.  Because $\L^2_{\w_{j+1}}\leq \L^2((i_1,...,\f_j^-))$, we have that $\L'^2_{\w_{j+1}}\leq \L^2((i_1,...,\f_{j+1}^-))$ by definition of fragment REFs and Property 3 of the network. So, by mathematical induction, a representation exists.

\subsection{Theorem 2}
This is stated as a proposition in \cite{varfixtheorem} though we could not find a proof. We prove it in alignment with the way Theorem 2 is stated.  Start by adding the constraint 
\begin{equation*}
x_j\geq 0\label{augeq}
\end{equation*}
to $P$ and assign the corresponding additional variable in the dual problem $\gamma$ the value of $\gamma^*=\bar{c}_j = c_j - y^\top E_j$.  The new primal problem is equivalent to the old and so is its linear relaxation.  Also, the new dual solution has the same cost because $\gamma$ has coefficient zero in the objective.  It is feasible because the constraint associated with $x_j$ becomes $y^\top E_j + \gamma \leq c_j$.  At $\gamma^* = c_j - y^\top E_j$, the left hand side is $y^\top E_j + c_j - y^\top E_j = c_j\leq c_j$.  Now change the right hand side of (\ref{augeq}) to one.  This increases the dual objective by $\gamma^* = \bar{c}_j$ and the dual solution remains feasible as no constraints in the dual change.  So the optimal dual solution has cost at least $z_{lb} + \bar{c}_j > z_{ub}$. Hence by duality, forcing $x_j = 1$ in $P$ increases the objective cost to at least $z_{lb}+\bar{c}_j>z_{ub}$.

\subsection{Theorem 3}
The result follows directly from the fact that $\rho_\w'$ is a lower bound on the reduced cost of any route with a representation containing a resourced copy of $\f_\w$.

\subsection{Theorem 4}
Let $\rr$ be the route with minimum reduced cost that has a representation containing $\w$.  Let that representation be $(\w_1,...,\w_k)$ and let the index of $\w$ in this sequence be $j$. 
Let $\hat{\p}_1 = \p_{{\w_1}}\oplus...\oplus \p_{{\w_{j-1}}}$ and $\hat{\p}_2 = \p_{{\w_{j+1}}}\oplus...\oplus \p_{{\w_k}}$ be the sub-paths before and after the resourced fragment's path in the route.  We consider the case where $\hat{\p}_1$ and $\hat{\p}_2$ are undominated in the forward and backward labelling algorithms respectively. Similar reasoning can be used when either or both are dominated by replacing $\hat{\p}_1$ and/or $\hat{\p}_2$ with the path that dominates them. Path 
$\hat{\p}_2\in \BB(\hat{\p}_1\oplus \p_{\w})$ otherwise $\rr$ would not be resource feasible.  Also, $\bL^2(\hat{\p}_2)\geq \tau_\w^2$ by the presupposition of the theorem. So
\begin{equation*}
    \bar{c}_{\hat{\p}_1} + \bar{c}_{\w} + \min_{ \p_2\in\BB(\hat{\p}_1\oplus \p_{\w})\text{ s.t. }\bL^2(\p_2) \geq \tau_\w^2}\bar{c}_{\p_2} \leq \bar{c}_{\hat{\p}_1} + \bar{c}_{\w} + \bar{c}_{\hat{\p}^2}= \rho_\w.\label{vf3ineq2}
\end{equation*}
We have that $\L^2(\hat{\p}_1)\in\Lambda_{\w^+}$ by the condition we place on representations in the theorem.  Also $\L^2(\hat{\p}_1)\in T_{\f_\w}$ because $\rr$ is resource feasible. This means $\hat{P}_1$ is one of the paths in the outer minimization of (16) and the result follows.

\section{DDD Algorithm}
Algorithm \ref{dddalg} describes how DDD is used to solve the VRP.  It considers the more complex case of two disposable resources; time $t$ and load $q$.  The simpler case of just including time $t$ is achieved by removing all mentions of load $q$ from the algorithm.  Line 1 initialises a flag that signifies whether an optimal solution has been found.  Then, while an optimal solution has not been found, lines 3 and 4 define and solve RFN to optimality.  If all chains in the solution are representations, we toggle the $Solved$ flag and terminate with an optimal solution.  Otherwise, we separate at least one minimal underestimating chain in line 6.  Algorithm \ref{dddalg} only collects minimal underestimating chains from the final RFN solution, however in practice, they can be collected from multiple/all solutions found. 
Lines 7 to 14 ensure the removal of each minimal underestimating chain. See Section 3.2 for details.  
\RestyleAlgo{ruled}
\LinesNumbered
\begin{algorithm}
\caption{Solve RFN}\label{dddalg}
\KwData{Resource Expanded Network $\N = (\V, \W\cup\A)$}
\KwResult{RFN solution $x_\w^*\in\{0,1\}$ for $\w\in\W$ corresponding to a set of representations}
$Solved\gets False$\;
\While{$Solved\neq True$}{
    Define RFN over $\N$\;
    Solve RFN by BB to obtain solution $x_\w^{*}$ for $\w\in\W$\;
    \eIf{\upshape There exists at least one chain that is not a representation}{
    Separate a non-empty set of minimal underestimating chains $\chi$\;
    \For{\upshape each chain $(\w_1,...,\w_k)\in\chi$}{
        Calculate $\V(\w_1,...,\w_k)$\;
        \For{$j=1,...,k-1$ and each unique load $q$ of a $\vv\in\V(\w_1,...,\w_k)$} {
            Add $\vv = ((\f_{\w_j}^-, (\tau^r_{\f_{\w_j}})_{r\in R_1}), (\tau^t_{\f_{\w_j}}\circ...\circ\tau^t_{\f_{\w_1}}(\alpha),\tau^q_{\f_{\w_j}}\circ...\circ\tau^q_{\f_{\w_1}}(q)))$ to $\V$\;
            Add all resourced fragments starting at $\vv$ to $\W$\;
            Update all arrival resourced nodes of resourced fragments that now violate the longest arc property\;
        }
    }
    }{
        $Solved\gets True$\;
    }
}
Return optimal solution $x_\w^*$ for $\w\in\W$\;
\end{algorithm} 

\section{PDPTW Implementation Details}\label{pdptwdetails}
\subsection{VRPSolver Parameters}
In step 1 of the overall algorithm we solve MP by CG using VRPSolver.  We do not use ng-neighbourhoods and so the initial solution can include all but 2-cycles.  If the initial $\lceil z_{lb}^v\rceil$ is the same as the number of vehicles in the initial heuristic upper bound solution, we finish.  Otherwise, we separate lm-R1Cs according to the VRPSolver parameter values given below.  We only include values that are different from those in \cite{vrpsolver}.

\begin{table}[!h]\centering
\caption{VRPSolver parameters for Step 1 of the overall algorithm}\label{step1vrpsolverparams}
\scriptsize
\begin{tabular}{lrr}\toprule
Parameter &Value \\\midrule
$\mathtt{RCSPstopCutGenTimeThresholdInPricing}$ &30 \\
$\mathtt{RCSPhardTimeThresholdInPricing}$ &60 \\
$\mathtt{RCSPmaxNumOfLabelsInEnumeration}$ &0 \\
$\mathtt{RCSPmaxNumOfEnumeratedSolutions}$ &0 \\
$\mathtt{RCSPmaxNumOfEnumSolutionsForMIP}$ &0 \\
$\mathtt{RCSPmaxNumOfEnumSolsForEndOfNodeMIP}$ &0 \\
$\mathtt{RCSPrankOneCutsMemoryType}$ &2 \\
$\mathtt{CutTailingOffCounterThreshold}$ &10 \\
$\mathtt{MaxNbOfBBtreeNodeTreated}$ &1 \\
$\mathtt{RCSPinitNGneighbourhoodSize}$ &0 \\
$\mathtt{RCSPmaxNGneighbourhoodSize}$ &0 \\
$\mathtt{RCSPapplyReducedCostFixing}$ &0 \\
$\mathtt{ReducedCostFixingThreshold}$ &0 \\
$\mathtt{SafeDualBoundScaleFactor}$ &1000000 \\
\bottomrule
\end{tabular}
\end{table}

The cuts are always separated when completing step 2 of the overall algorithm.  The only parameter value that is changed is $\mathtt{CutTailingOffCounterThreshold}=1$.

\subsection{Details of Solving the Fragment-Based Formulation}
The forward and backward labels are not made public in BaPCod, so we implemented our own slower labelling algorithms to backwards-engineer these from the final dual solution of MP found by VRPSolver.  Our implementation was slower than BaPCod's, especially for the MP dual solutions found in step 1 of the overall algorithm that required lm-R1Cs: those for CC70 and DD60 where the initial heuristic upper bound solution contained more vehicles than $\lceil z_{lb}^v\rceil$.  For these two instances, we saved the dual solution found by BaPCod with no lm-R1Cs, and then every five iterations of cut separation after this.  We then performed formulation leveraging with these dual solutions consecutively which reduced the overall time taken.  Each iteration of solving RFN by BB in DDD was terminated whenever a solution with underestimating chains and cost less than the best known solution was found.  This was because initially building up the network quickly was found to be important.  Resourced nodes were added for every minimal underestimating chain of minimal length found in every solution found.  When performing CERE, we chose $\mathtt{MaxIncrease'}=500$ and $\mathtt{MaxIncrease''} = 6000$.

\section{TDDRP Labelling Algorithm}\label{tddrplabelling}
This section describes the details of the forward labelling algorithm used to solve PP for the TDDRP. The labelling deviates significantly from the classical labelling described in the main paper.  This is because of the synchronisation complexity between the truck and the drone.  
We define the REFs for $t_\text{truck}$, $t_\text{drone}$ and $q$.  When the truck and drone are at the same location, we extend the label by each drone flight $\f^\text{drone}\in\F^\text{drone}_{i_\text{drone}}$ with $\max\{t_\text{drone}, t_\text{truck}\}\leq\mathtt{latest}_{\f^\text{drone}}$:
\begin{equation}
f_{\f^\text{drone}}^{t_\text{drone}}(t_\text{truck}, t_\text{drone}) = \mathtt{finish}_{\f^\text{drone}}(\max\{t_\text{drone},t_\text{truck}\}).
\end{equation}
The truck can deploy the drone at a location before or after the request's time window so when a label is extended by the truck moving to $j\in (N\cup \{i_\text{drone}\})\setminus\{i_\text{truck}\}$, we ensure $t_\text{truck}$ represents the truck's earliest arrival at the new location:
\begin{equation}
    f_{i_\text{truck}j}^{t_\text{truck}}(i_\text{truck}, i_\text{drone}, t_\text{truck}, t_\text{drone}) = \begin{cases}
        \max\{t_\text{truck},t_\text{drone},a_{i_\text{truck}}\} + t^\text{truck}_{i_\text{truck}j} &\text{ if }i_\text{truck}=i_\text{drone},\\
        \max\{t_\text{truck},a_{i_\text{truck}}\} + t^\text{truck}_{i_\text{truck}j}&\text{ otherwise}.
    \end{cases}
\end{equation}
The used truck capacity is also updated each time the truck moves or the drone is deployed:
\begin{equation}
    f_{i_\text{truck}j}^q(q) = q + d_j\text{ and }f_{\f^\text{drone}}^q(q) = q+d_{\f^\text{drone}}.
\end{equation}

The algorithm initialises with label $(0,0,0,(0,0,0))$.  The feasible extensions of label $\l = (i_\text{truck}, i_\text{drone}, \bc, (t_\text{truck}, t_\text{drone}, q))$ are described below by considering the cases where $i_\text{truck}=i_\text{drone}$ and $i_\text{truck}\neq i_\text{drone}$.
\subsection{Case 1: $i_\text{truck}=i_\text{drone}\in N\cup\{0\}$}  
\begin{enumerate}
    \item The truck can carry the drone to a new request $j\in (N\cup\{n+1\})\setminus\{i_\text{truck},i_\text{drone}\}$ in which case the new label is

\begin{equation}
    \l' = (j, j, \bc', (f_{i_\text{truck}j}^{t_\text{truck}}(i_\text{truck},i_\text{drone},t_\text{truck},t_\text{drone}),\\ 
    f_{i_\text{drone}j}^{t_\text{truck}}(i_\text{truck},i_\text{drone},t_\text{truck},t_\text{drone}), f_{i_\text{truck}j}^q(q))),\label{carry}
\end{equation}
where the updated reduced cost is $\bc'=\bc + \pi_{i_\text{truck}}/2 + \pi_j/2$.  This is resource feasible if $f_{i_\text{truck}j}^{t_\text{truck}}(i_\text{truck},i_\text{drone},t_\text{truck},t_\text{drone})\leq b_{j}$ and $f_{i_\text{truck}j}^q(q)\leq Q^\text{truck}$.
\item The drone can be deployed on drone flight $\f^\text{drone}=(i_1,...,i_k)\in\F^\text{drone}_{i_\text{drone}}$ satisfying $\max\{t_\text{truck},t_\text{drone}\}\leq\mathtt{latest}_{\f^\text{drone}}$.  This extension yields 
\begin{equation}
    \l' = (i_\text{truck}, i_k, \bc-\sum_{j=2}^{k-1}\pi_{i_j},(t_\text{truck}, f_{\f^\text{drone}}^{t_\text{drone}}(t_\text{truck},t_\text{drone}), f_{\f^\text{drone}}^q(q))).\label{separate}
\end{equation}
This is resource feasible if  $f_{\f^\text{drone}}^q(q) \leq Q^\text{truck}$. 
\end{enumerate}

\subsection{Case 2: \textbf{$i_\text{truck}\neq i_\text{drone}$}}

In this case only the truck can move to a new location $j\in (N\cup \{i_\text{drone}\})\setminus\{i_\text{truck}\}$.  This extension yields new label
\begin{equation}
    \l' = (j, i_\text{drone}, \bc', f_{i_\text{truck}j}^{t_\text{truck}}(i_\text{truck},i_\text{drone},t_\text{truck},t_\text{drone}), t_\text{drone}, f_{i_\text{truck}j}^q(q))\label{mvtk}
\end{equation}
and has the same resource feasibility conditions as (\ref{carry}).

The backward labelling algorithm is similar to the forward labelling algorithm, however, $i_\text{truck}$ and $i_\text{drone}$ become different by the truck moving to another location, rather than the drone being deployed on a drone flight. Hence, a drone flight is then used to synchronise the pair at a subsequent truck location. 

\section{TDDRP Fragment Enumeration}\label{tddrpfragenum}
Type 1 fragments can be fully enumerated as there are at most $n^2$ of these.  On the other hand, type 2 fragments cannot be enumerated for large instances.  We consequently perform FL on partial type 2 fragments with a full drone flight and partial vehicle leg $\f = (\f^\text{drone} = (i_1^1,...,i_{k_1}^1),\f^\text{truck}=(i_1^2 = i_1^1,...,i_l^2))$.  Since forward and backward labelling are significantly more complex in the TDDRP, and $i^2_l \neq i_{k_1}^1$ unless $\f^\text{truck}$ is a full truck leg, we introduce a new backward label set to help explain. Let $\BB(i^\text{truck}, i^\text{drone}, t^\text{truck}, t^\text{drone}, q)$ be the set of undominated backward labels where the truck and drone start at $i^\text{truck}$ and $i^\text{drone}$ respectively with resource consumption values component-wise greater than or equal to $(t^\text{truck}, t^\text{drone}, q)$, and arrive at $n+1$ feasibly. We only consider forward labels where the truck and drone are at the same location because $i_1^1=i_1^2$.  So let $\FF(i,\Lambda)$ be the set of forward labels where the truck and drone leave $0$ and both arrive at $i$ by minimum time $t$ with cumulative load $q$ where $(t,q)\in \Lambda$.

We now describe how we calculate a lower bound $\rho$ on the reduced cost of routes that contain partial fragment $\f=(\f^\text{drone}, \f^\text{truck})$.  
The reduced cost contribution of $\f$ is 
\begin{equation*}
    \bar{c}_\f = \sum_{j=1}^{l-1}c_{i_j^2 i_{j+1}^2}^\text{truck}+\sum_{j=1}^{k_1-1}c_{i_j^2 i_{j+1}^2}^\text{drone}-\sum_{j=2}^{k_1-1}\pi_{i^1_j}-\sum_{j=1}^{l}\pi_{i^2_j}/2.
\end{equation*}
For each label $\l\in\FF(i_1^1,T_\f)$, we perform the drone deployment extension (\ref{separate}) to calculate $t^\text{drone}$ and $q$ given the drone completes drone flight $\f^\text{drone}$ after the partial route represented by $\l$.  
We perform repeated vehicle movement extensions (\ref{mvtk}) to calculate $t^\text{truck}$ and $q$ given the truck then moves from $i_1^2$ to $i_l^2$.  
Then for each backward label $\bl\in \BB(i_l^2,i_{k_1}^1,t^\text{truck}, t^\text{drone},q)$, $\rho$ is updated to
\begin{equation}
    \rho = \min\{\rho, \bar{c}_{\l} + \bar{c}_\f + \bar{c}_{\bl}\},
\end{equation}
where $\bc_\l$ and $\bc_{\bl}$ are the reduced costs of forward label $\l$ and backward label $\bl$ respectively.  The final value of $\rho$ after all forward paths in $\FF(i_1^1,T_\f)$ are searched is a lower bound on the reduced cost of routes that contain partial fragment $\f$.  If $\rho$ is greater than the optimality gap, $\f$ is discarded from our fragment enumeration.

\bibliographystyle{supplemental_material}
\bibliography{supplemental_material}